\documentclass[11pt,twoside]{amsart}
\usepackage{amssymb}

\usepackage[latin1]{inputenc}
\usepackage[T1]{fontenc}
\usepackage{mathtools}
\usepackage{graphicx}
\usepackage{tikz}
\usepackage{mathabx}
\usetikzlibrary{chains}
\usepackage[OT2,T1]{fontenc}
\usepackage{enumitem}
\PassOptionsToPackage{pdfusetitle,pagebackref,colorlinks}{hyperref}
\usepackage{bookmark}
\hypersetup{
  linkcolor={red!70!black},
  citecolor={green!70!black},
  urlcolor={blue!80!black}
}

\usepackage[latin1]{inputenc}
\usepackage{amsmath}
\usepackage{amsthm}

\usepackage[all]{xy}
\usepackage{hyperref}
\newtheorem{thm}{Theorem}[section]
\newtheorem{thmx}{Theorem}

\newtheorem{prop}[thm]{Proposition}

\newtheorem{lem}[thm]{Lemma}
\newtheorem{cor}[thm]{Corollary}
\newtheorem{conj}[thm]{Conjecture}
\numberwithin{equation}{section}

\theoremstyle{definition}

\newtheorem{remark}[thm]{Remark}

\DeclareFontFamily{U}{mathc}{}
\DeclareFontShape{U}{mathc}{m}{it}%
{<->s*[1.03] mathc10}{}
\DeclareMathAlphabet{\mathcal}{U}{mathc}{m}{it}
\newcommand{\kend}{\mathcal{E\mkern-3mu nd}}

\newcommand{\Db}{{\rm D}^{\rm b}}
\newcommand{\ch}{{\rm ch}}

\newcommand{\Br}{{\rm Br}}

\newcommand{\Coh}{{\rm Coh}}
\newcommand{\per}{{\rm per}}
\newcommand{\ind}{{\rm ind}}

\newcommand{\NS}{{\rm NS}}
\newcommand{\Pic}{{\rm Pic}}

\newcommand{\rk}{{\rm rk}}

\newcommand{\Sym}{{ S}}

\newcommand{\cal}{\mathcal}

\newcommand{\kc}{{\cal C}}

\newcommand{\km}{{\cal M}}
\newcommand{\ko}{{\cal O}}
\newcommand{\kp}{{\cal P}}

\newcommand{\kx}{{\cal X}}

\newcommand{\ZZ}{\mathbb{Z}}
\newcommand{\QQ}{\mathbb{Q}}
\newcommand{\RR}{\mathbb{R}}
\newcommand{\CC}{\mathbb{C}}

\newcommand{\PP}{\mathbb{P}}

\DeclareSymbolFont{cyrletters}{OT2}{wncyr}{m}{n}
\DeclareMathSymbol{\Sha}{\mathalpha}{cyrletters}{"58}

\renewcommand{\to}{\xymatrix@1@=15pt{\ar[r]&}}
\newcommand{\lto}{\xymatrix@1@=15pt{&\ar[l]}}
\renewcommand{\rightarrow}{\xymatrix@1@=15pt{\ar[r]&}}
\renewcommand{\mapsto}{\xymatrix@1@=15pt{\ar@{|->}[r]&}}
\newcommand{\mapslto}{\xymatrix@1@=15pt{&\ar@{|->}[l]&}}
\renewcommand{\twoheadrightarrow}{\xymatrix@1@=18pt{\ar@{->>}[r]&}}

\renewcommand{\hookrightarrow}{\xymatrix@1@=15pt{\ar@{^(->}[r]&}}
\newcommand{\hook}{\xymatrix@1@=15pt{\ar@{^(->}[r]&}}

\newcommand{\congpf}{\xymatrix@1@=15pt{\ar[r]^-\sim&}}
\renewcommand{\cong}{\simeq}

\newcommand{\EV}[1]{}

\newcommand{\Old}[1]{}
\makeatletter
\def\blfootnote{\xdef\@thefnmark{}\@footnotetext}
\makeatother

\begin{document}

\title[]{The period--index problem for  hyperk\"ahler varieties: Lower and upper bounds}

\author[A.\ Bottini \& D.\ Huybrechts]{Alessio Bottini \& Daniel Huybrechts}

\address{Mathematical Institute and Hausdorff Center for Mathematics,
University of Bonn, Endenicher Allee 60, 53115 Bonn, Germany}
\email{bottini@math.uni-bonn.de \& huybrech@math.uni-bonn.de}

\begin{abstract}  \vspace{-2mm} It is expected that a stronger form of the period-index conjecture
holds for hyperk\"ahler varieties. Following ideas of Hotchkiss, we provide further evidence  for this expectation by proving a version in which the index is replaced by the Hodge-theoretic index. 

We also show that the hyperk\"ahler period-index conjecture is optimal. As an application, we prove that Mumford--Tate general  hyperk\"ahler varieties cannot be covered by families of elliptic curves passing through a fixed point. 

By extending work of Hotchkiss, Maulik, Shen, Yin, and Zhang, we prove
the hyperk\"ahler period-index conjecture for non-special coprime Brauer class  on
hyperk\"ahler varieties of ${\rm K3}^{[n]}$-type without any restriction on the Picard number.
\end{abstract}

\maketitle
\blfootnote{Both authors have been supported by the ERC Synergy Grant HyperK (ID 854361).}

\section{Introduction}
Before turning to the period-index conjecture for hyperk\"ahler varieties and outlining the results proved in this paper, we briefly review the general conjecture.

\subsection{Period-index} With any Brauer class $\alpha\in\Br(X)$ on a smooth projective variety $X$ one associates
two numerical invariants, its period $\per(\alpha)$ and its index $\ind(\alpha)$.
The period is simply the order of $\alpha$ as an element of the abelian group $\Br(X)$ while
the index is the minimal $m$ such that the restriction of $\alpha$ to some non-empty Zariski
open subset $U\subset X$ can be represented by a Brauer--Severi variety $P\to U$ of relative dimension $m-1$. According to a classical result  going back to Albert \cite{Albert}, one knows that $\per(\alpha)\mid\ind(\alpha)\mid\per(\alpha)^{e(\alpha)}$ for some integer $e(\alpha)$ which a priori depends on the specific class
$\alpha$. \smallskip

It is generally believed that the exponent $e(\alpha)$ can be chosen to only depend on $X$. This is made
more precise by the following, cf.\ \cite{CT}:
 
\begin{conj}[Colliot-Th\'el\`ene]\label{conj:PI}
All Brauer classes $\alpha\in\Br(X)$ on a smooth projective variety $X$ over an algebraically closed field
satisfy 
\begin{equation}\label{eqn:PI}
\ind(\alpha)\mid\per(\alpha)^{\dim(X)-1}.
\end{equation}
\end{conj}

The conjecture is wide open in general, but a few things are known. For example, the `discriminant avoidance' result of de Jong and
Starr \cite{dJS} shows that once Conjecture \ref{conj:PI} is established the upper bound (\ref{eqn:PI})
also holds for all `ramified' Brauer classes, i.e.\ classes $\alpha\in\Br(K(X))$ over the function field, cf.\ the recent more geometric account \cite{GH}. Also, due to joint work with Mattei \cite{HM}, it is known that a  uniform upper bound exists. More precisely, there exists an integer $e(X)$ only depending on $X$ with $\ind(\alpha)\mid\per(\alpha)^{e(X)}$
for all `unramified' $\alpha\in\Br(X)$. A similar result was established by de Jong and Perry \cite{dJP} conditionally on the standard Lefschetz conjecture in degree two. At the moment, 
lacking a sufficiently strong discriminant avoidance result, the existing uniform bounds do not extend to ramified classes. Famously, (\ref{eqn:PI}) is known for curves by Tsen's theorem and
for surfaces by a result of de Jong, see \cite{dJSurf,dJS,L2}. It is also known that if true the conjecture is optimal. Examples proving optimality, especially of unramified classes, are not easy to come by. The first one is due to Gabber, see the appendix to \cite{CT3}; it has  recently been generalised in \cite{dJP,HotchHodge}.\smallskip

The best evidence for Conjecture \ref{conj:PI} and the first conceptual (geometric) explanation of the 
particular
shape of the upper bound is provided by recent results of Hotchkiss \cite{HotchHodge}, see
also \cite{dJP}. He defines the Hodge-theoretic index $\ind_H(\alpha)$ for any unramified Brauer class $\alpha\in\Br(X)$, see \S\! \ref{sec:HoHo}, and shows 
\begin{equation}\label{eqn:PIH}
\ind_H(\alpha)\mid\per(\alpha)^{\dim(X)-1}
\end{equation} for all
classes $\alpha\in\Br(X)$ for which $\per(\alpha)$ and $(\dim(X)-1)!$ are coprime. Moreover,
he argues that if the integral Hodge conjecture holds for a Brauer--Severi variety $P\to X$ representing $\alpha$, then also the original (\ref{eqn:PI}) holds. Indeed, in this case the general
$\ind_H(\alpha)\mid\ind(\alpha)$ becomes an equality. A similar strategy has been applied
by de Jong and Perry \cite{dJP}.

\begin{remark}
In a number of recent papers \cite{HotchHodge,HMSYZ,HuyPI,HM}, certain results are not established for all classes $\alpha\in\Br(X)$
but only for those with $\per(\alpha)$  coprime to a fixed positive integer $N_X$.  For short, we say that in this case $\alpha$ is coprime to $N_X$.  The integer $N_X$ only depends on $X$
and often only on certain numerical invariants of $X$,
e.g.\ in the result of Hotchkiss one chooses $N_X=(\dim(X)-1)!$.
The exact value of $N_X$ is usually not important and rarely optimal. In order to avoid 
having to specify $N_X$,  we often simply say  that a certain result holds for all \emph{coprime classes} $\alpha\in\Br(X)$.
\end{remark}


\subsection{Period-index for hyperk\"ahler varieties}

Although the upper bound (\ref{eqn:PI}) is known to be optimal in general, it is conceivable that better bounds
exist for particular classes of varieties. Of course, this is only relevant if the Brauer group is not trivial altogether. One particular class of varieties with an interesting Brauer group is provided by K3 surfaces and  higher-dimensional generalisations thereof, hyperk\"ahler varieties. In this paper, a \emph{hyperk\"ahler variety} will always refer to a simply connected, smooth, projective, complex variety $X$ with $H^{2,0}(X) $ spanned by an
everywhere non-degenerate two-form. For these, the following stronger version of Conjecture \ref{conj:PI} was put forward in \cite{HuyPI}.

\begin{conj}\label{conj:PIHK}
All Brauer classes $\alpha\in\Br(X)$ on a hyperk\"ahler variety $X$ satisfy 
\begin{equation}\label{eqn:PIHK}
\ind(\alpha)\mid\per(\alpha)^{\dim(X)/2}.
\end{equation}
\end{conj}

The conjecture is known to hold for all hyperk\"ahler varieties in various dense countable unions of 
codimension one loci in the moduli space of hyperk\"ahler varieties. We summarise the previously established cases:\smallskip

\begin{enumerate}[topsep=0pt,itemsep=0.6ex]
\item[$\bullet$] If $X$ admits a Lagrangian fibration $X\to\PP^n$, then (\ref{eqn:PIHK})
holds for all coprime classes $\alpha\in \Br(X)$, see \cite[Thm.\ 0.3]{HuyPI}.
\item[$\bullet$] If $X$ is the Hilbert scheme $S^{[n]}$ of a K3 surface $S$, then (\ref{eqn:PIHK})
holds for all classes $\alpha\in \Br(X)$, see \cite[Thm.\ 0.4]{HuyPI}.
\item[$\bullet$] If $X$ is a hyperk\"ahler variety of ${\rm K3}^{[n]}$-type of Picard number
at least two,  then (\ref{eqn:PIHK}) holds for all non-special coprime classes $\alpha\in \Br(X)$, see \cite[Thm.\ 0.5]{HMSYZ}.
\end{enumerate}\smallskip

Hotchkiss et al \cite{HMSYZ} also establish the uniform bound $e(X)=\dim(X)$ for coprime Brauer classes on arbitrary hyperk\"ahler varieties of ${\rm K3}^{[n]}$-type. The additional assumption of $\alpha$ being non-special is a condition on the square of the $B$-field representing $\alpha$,
see Remark \ref{rem:ellper}. \smallskip

Further examples providing evidence for Conjecture \ref{conj:PIHK} can be obtained
by considering surface decomposable hyperk\"ahler varieties, cf.\ \cite[\S\! 3.5]{HuyPI}.
According to results of Voisin, see \cite[Thm.\ 3.3]{VoisinDecom}, the Fano variety of lines on a cubic fourfold, which is of ${\rm K3}^{[2]}$-type, has a surface decomposition. Hence, 
(\ref{eqn:PIHK}) holds for coprime classes on Fano varieties of lines, but see Theorem \ref{thm:thm4} below for a more general result. More interestingly, also the LSV compactification of the intermediate Jacobian
fibration associated with a cubic fourfold \cite{LSV}, a  hyperk\"ahler variety of OG10-type, is surface decomposable and, therefore,
(\ref{eqn:PIHK}) holds for all coprime classes. There are also examples of surface
decomposable hyperk\"ahler varieties of OG6-type, see \cite{MRS}.
\smallskip

In the same sense that Hotchkiss's result $\ind_H(\alpha)\mid\per(\alpha)^{\dim(X)-1}$  for an arbitrary smooth projective variety $X$ is
the strongest evidence yet for the period-index conjecture (\ref{eqn:PI}), we view our first theorem  as convincing
evidence for the hyperk\"ahler version (\ref{eqn:PIHK}). The proof will be given in \S\! \ref{sec:ProofThm1}.

\begin{thmx}\label{thm:1}
For all coprime classes $\alpha\in\Br(X)$ on a hyperk\"ahler variety $X$ 
one has $$\ind_H(\alpha)\mid\per(\alpha)^{\dim(X)/2}.$$
In addition, if $\alpha$ is represented by a Brauer--Severi variety $P\to X$ and the integral Hodge conjecture holds for $P$, then Conjecture \ref{conj:PIHK} holds.
\end{thmx}

In fact, in \S\! \ref{sec:ProofThm1} we will prove that all Hodge-theoretic indices defined in terms of the various Hodge structures 
$H^\ast$, $SH$, and $\Sym^n\widetilde H$ naturally associated with a hyperk\"ahler
variety, see \S\! \ref{sec:HIIntro} below, satisfy the same divisibility.

\subsection{Optimal period-index}

The next result shows that Conjecture \ref{conj:PIHK} is optimal. More is true, we cannot
only exclude a stronger upper bound to hold for all Brauer classes, but we show in fact that it cannot even exist 
when one restricts to coprime classes (with respect to any integer $N_X$). The next result will be proven in \S\! \ref{sec:PrThmB}.

\begin{thmx}\label{thm:2}
Assume $X$ is a general hyperk\"ahler variety of Picard rank $\rho(X)\leq b_2(X)-6$. Then for any non-special coprime class
$\alpha\in \Br(X)$ one has
\begin{equation}\label{eqn:thm2}
\per(\alpha)^{\dim(X)/2}\mid \ind(\alpha).
\end{equation}
\end{thmx}
In \S\! \ref{sec:Lower} we will make the assumption more precise; `general' is specified as `Mumford--Tate general'.  Note that for the known deformation types both assumptions
hold away from a countable union of proper closed subsets in the moduli space.
A fairly complete picture emerges in the case of hyperk\"ahler varieties  admitting a Lagrangian fibration and hyperk\"ahler varieties of ${\rm K3}^{[n]}$-type, cf.\ Corollary \ref{cor:equality} and Corollary \ref{cor:LagMT}.
\subsection{Hodge-theoretic indices}\label{sec:HIIntro}
The key to the above two theorems is the additional symmetry of the Hodge structure of
a hyperk\"ahler variety which gives rise to refinements of Hotchkiss's Hodge index $\ind_H(\alpha)$.
The Mukai lattice $\widetilde H$, its symmetric power $\Sym^n\widetilde H$, and
the Verbitsky component $SH(X)\subset H^\ast (X)$, i.e.\ the sub-algebra of the cohomology generated by
$H^2$, all lead to numerical invariants
$\ind_{\widetilde H}(\alpha)$, $\ind_{\Sym^n\widetilde H}(\alpha)$, $\ind_{H^\ast}(\alpha)$, and $\ind_{SH}(\alpha)$
of any topologically trivial Brauer class $\alpha\in\Br(X)$, see \S\! \ref{sec:HI} for precise definitions.
They are linked to each other as follows, see \S \! \ref{sec:HI} for the verifications.

\begin{thmx}\label{thm:thm3}
Let $X$ be a hyperk\"ahler variety of dimension $2n$. Then $\ind_H(\alpha)\mid\ind(\alpha)$,
\begin{eqnarray*}
\ind_{H^\ast}(\alpha)\mid\gcd\left(
(2n)!\cdot \ind_H(\alpha), \ind_{SH}(\alpha)\right),\\~\text{\rm and }~
\ind_{\Sym^n\widetilde H}(\alpha)\mid\text{\rm gcd}\left(
n!\cdot\ind_{SH}(\alpha),\ind_{\widetilde H}(\alpha)^n\right).
\end{eqnarray*}
For coprime Brauer classes $\alpha\in\Br(X)$ one has
$$\ind_H(\alpha)=\ind_{H^\ast}(\alpha)=\ind_{SH}(\alpha)=\ind_{\Sym^n\widetilde H}(\alpha)\mid\ind_{\widetilde H}(\alpha)^n.$$
Moreover, $\per(\alpha)$ divides all four indices.
\end{thmx}

\subsection{Prime Brauer classes on ${\rm K3}^{[n]}$} Combining the arguments in \cite{HMSYZ}, which
crucially rely on results of Markman \cite{Markman}, with work of O'Grady \cite{OG2}, we give further evidence for the period-index conjecture by the following result.

\begin{thmx}\label{thm:thm4} Assume $X$ is a hyperk\"ahler variety of ${\rm K3}^{[n]}$-type.
If $\alpha\in \Br(X)$ is a non-special coprime Brauer class, then
\begin{equation}\label{eqn:thmD}
\ind(\alpha)\mid\per(\alpha)^{\dim(X)/2}.
\end{equation}
For a special coprime Brauer class $\alpha$ {\rm (\ref{eqn:thmD})} holds true
if its period $\per(\alpha)$ has a reduced prime factor decomposition  $\per(\alpha)=\prod p_i$, $p_i\ne p_j$. 
\end{thmx}

The theorem extends \cite[Thm.\ 0.5]{HMSYZ} from the case of Picard number at least two to all hyperk\"ahler varieties of ${\rm K3}^{[n]}$-type. Its proof can be found in \S\S\! \ref{sec:special} and \ref{sec:nonspecial}.\smallskip

We also note the following combination of a version of Theorem \ref{thm:2},
see Remark \ref{rem:K3n}, and Theorem \ref{thm:thm4}.

\begin{cor}\label{cor:equality}
Assume $X$ is a Mumford--Tate general hyperk\"ahler variety of ${\rm K3}^{[n]}$-type of Picard number one.  Then
$$\per(\alpha)^{n}=\ind(\alpha)$$
and $\alpha\in\Br(X)$ is a non-special coprime Brauer class.\qed
\end{cor}

\subsection{Covering families}
In \cite[Question 1.2]{Voisin} Voisin asks whether a very general hyperk\"ahler variety is covered
by elliptic curves. The motivation comes from recent studies of birational invariants of projective varieties like covering gonality and degree of irrationality. Combining the results in \cite{HuyPI,HM} with Theorem \ref{thm:2} we show that the period-index problem affects this question, see \S\! \ref{sec:proofthm5} for the proof.

\begin{thmx}\label{thm:thm5}
Assume $X$ is a Mumford--Tate general hyperk\"ahler variety and $\{\kc_t\}$ is a dominating family
of curves in $X$ all passing through one fixed point $x\in X$. Then the geometric genus satisfies
$$g(\kc_t)\geq \dim(X)/2.$$
\end{thmx}

In particular, general hyperk\"ahler varieties of dimension at least four do not admit covering families
of elliptic curves  all passing through one fixed point.

\medskip

\noindent
{\bf Acknowledgement:} We wish to thank  Moritz Hartlieb for a stimulating discussion related to 
\S\! \ref{sec:curvesHilb}. This work started while the first named author was funded by the 
 Max Planck Institute for Mathematics, Bonn.
The first version of the paper was written while the second named author enjoyed the hospitality of Imperial College, London. 

\section{Hodge-theoretic indices}\label{sec:HI} With a hyperk\"ahler variety $X$ various Hodge structures are naturally associated. We will concentrate on the following ones (all to be recalled below):
\begin{equation}\label{eqn:HS}H^{2\ast}(X,\ZZ),~SH(X,\ZZ),~\widetilde H(X,\ZZ),
\text{ and } ~\Sym^n \widetilde H(X,\ZZ).
\end{equation}
Here, $SH(X,\ZZ)$ is the Verbitsky component and $\widetilde H(X,\ZZ)=H^2(X,\ZZ)\oplus U$ is the Hodge--Mukai lattice.
Using Tate twists, they will all be considered as Hodge structures of weight zero and torsion
will be ignored throughout. There are variants of the above where instead of integral cohomology
one uses topological $K$-theory. This is only really relevant for the first two. \smallskip 

The underlying idea to define the index of a Brauer class $\alpha\in\Br(X)$ with respect to any of these Hodge structures is always the same:  \smallskip

$\bullet$ Firstly, one uses the (topologically trivial) Brauer class $\alpha\in \Br(X)$ to `twist' any of the four integral Hodge structures $V$ in (\ref{eqn:HS}) to obtain a new integral Hodge structure $V_\alpha$
such that $V$ and $V_\alpha$ are commensurable, i.e.\ $V\otimes\QQ\cong V_\alpha\otimes\QQ$ as Hodge structures.\smallskip

$\bullet$ Secondly, for any integral Hodge structure $V$ of weight zero together with a given rank function $\rk\colon V\to \ZZ$ one defines $\ind(V)$ as the positive generator of the image of $\rk$ restricted to the subgroup of Hodge classes  $V_{\rm Hdg}\coloneqq V\cap V^{0,0}\subset V$, i.e.\ $\rk(V_{\rm Hdg})=\ind(V)\cdot \ZZ$.
The dependence of $\ind(V)$ on the rank function is suppressed in the notation.\smallskip

$\bullet$ Thirdly, if an integral Hodge structure $V$ comes with a rank function, then
 a natural rank function on the twisted Hodge structure $V_\alpha$ is given via the isomorphism
$V\otimes\QQ\cong V_\alpha\otimes\QQ$. One then defines the (a priori rational but often integral) number
$$\ind_V(\alpha)\coloneqq\ind(V_\alpha).$$

\subsection{Twisting the Hodge--Mukai lattice}\label{sec:twistHM}
We consider $H^2(X,\ZZ)$ of a hyperk\"ahler variety $X$
together with its  Hodge structure of weight two and the BBF-form $q$. However,  the following is purely algebraic and could be set up without a hyperk\"ahler variety in the background.\smallskip

The Hodge--Mukai lattice is  defined as the orthogonal sum $$\widetilde H(X,\ZZ)\coloneqq H^2(X,\ZZ)(1)\oplus U.$$
Here, $U$ is the hyperbolic plane, for which we fix a basis $e,f$ of isotropic vectors with $(e.f)=-1$
and for $b\in H^2(X,\ZZ)$ we often write $(r,b,s)$ instead of $r\,e+b+s\,f$.  The quadratic form on the Mukai lattice shall also be denoted $q$. The weight zero Hodge structure  of
 the Tate twist $H^2(X,\ZZ)(1)$ is extended  by declaring $U$ to be of Hodge type $(0,0)$.
 A natural rank function $\rk\colon \widetilde H(X,\ZZ)\to\ZZ$ is given by the projection
 to $\ZZ e$ or, alternatively, as $\rk=-q(~,f)$.\smallskip

Next, to a fixed topologically trivial Brauer class $\alpha\in\Br(X)$, we associate a twist 
$\widetilde H(X,\alpha,\ZZ)$ of the Hodge--Mukai lattice. For this, one has to fix a $B$-field lift of $\alpha$, i.e.\ a class $B\in H^2(X,\QQ)$ that under the exponential map is mapped to
$\alpha\in\Br(X)\subset H^2(X,\ko_X^\ast)$. One then defines $\widetilde H(X,B,\ZZ)$ as the Mukai lattice
with the Hodge structure determined by the following conditions: 
$H^{1,-1}(X,B)$ is spanned by $(0,\sigma,-q(\sigma,B))=\sigma-q(\sigma,B) \,f\in H^2(X,\CC)\oplus \CC f$ and $H^{0,0}(X,B)$ is the orthogonal complement of
$H^{1,-1}(X,B)\oplus H^{-1,1}(X,B)$. Here,  $\sigma\in H^{2,0}(X)$ is an arbitrary generator. It is known that the isomorphism type
of $\widetilde H(X,B,\ZZ)$ only depends on $\alpha$, see \cite{HS} for more details and motivation.\smallskip

As the lattice underlying $\widetilde H(X,B,\ZZ)$ is just the untwisted $\widetilde H(X,\ZZ)$,
we can keep the rank function $\rk\colon \widetilde H(X,B,\ZZ)\to \ZZ$, i.e.\ the projection to $\ZZ e$.
The same rank function can also be obtained by composing the original one with the
isomorphism of rational Hodge structures $\exp(B)\colon\widetilde H(X,B,\QQ)\cong \widetilde H(X,\QQ)$ 
given by multiplication with $\exp(B)=(1,B,B^2/2)$.
Following the general recipe we now define
$$\ind_{\widetilde H}(\alpha)\coloneqq \ind(\widetilde H(X,B,\ZZ)).$$
The definition is independent of the choice of the $B$-field lift $B$.
\subsection{Symmetric power of the Hodge--Mukai lattice}\label{sec:SymHM}
Let $X$ be a hyperk\"ahler variety
and $n$ a positive integer. We will eventually only be interested in the case that $2n=\dim(X)$.
 Then the symmetric power
$\Sym^n\widetilde H(X,\ZZ)$ is naturally endowed with a Hodge structure of weight zero and
also with a quadratic form, called $q^{(n)}$, with the property that
$q^{(n)}(a^{(n)},b^{(n)})=n!\,q(a,b)^n$ for the
symmetric powers $a^{(n)},b^{(n)}$ of any two classes $a,b\in \widetilde H(X,\ZZ)$. We define the rank function $\rk\colon \Sym^n\widetilde H(X,\ZZ)\to \ZZ$ as $\rk(v)=(-1)^nq^{(n)}(v,f^{(n)})$.
\smallskip 

There are a priori two ways to twist $\Sym^n\widetilde H(X,\ZZ)$ with respect to a topologically trivial Brauer class $\alpha\in\Br(X)$, for which we again fix a $B$-field lift $B\in H^2(X,\QQ)$. First, one could use the identification of lattices $\widetilde H(X,\ZZ)=\widetilde H(X,B,\ZZ)$ and the induced isomorphism between their symmetric powers, to actually define $\Sym^n\widetilde H(X,B,\ZZ)$ as the symmetric power of $\widetilde H(X,B,\ZZ)$. Second, $\Sym^n\exp(B)\colon\Sym^n\widetilde H(X,\QQ)\congpf \Sym^n\widetilde H(X,\QQ)$ is an isomorphism of vector spaces and we declare $\Sym^n\widetilde H(X,B,\ZZ)$ to be
the Hodge structure on $\Sym^n\widetilde H(X,\ZZ)\subset\Sym^n\widetilde H(X,\QQ)$ obtained by pull-back. It is easy to see that
these two definitions coincide.\smallskip

With either definition, we introduce yet another Hodge theoretic index $$\ind_{\Sym^n\widetilde H}(\alpha)\coloneqq\ind(\Sym^n\widetilde H(X,B,\ZZ)),$$
which is once more independent of the choice of $B$.

\subsection{Verbitsky component}\label{sec:VerbComp}
Let $X$ be a hyperk\"ahler variety of dimension $2n$. Consider the integral subalgebra
$SH(X,\ZZ)\subset \bigoplus_kH^{2k}(X,\ZZ)(k)$ generated by 
$H^2(X,\ZZ)(1)$ with its natural weight zero Hodge structure. By a result of Verbitsky \cite{Verb}, we
know that the natural map $S^\ast H^2(X,\ZZ)(1)\twoheadrightarrow
SH(X,\ZZ)$ is an isomorphism in degree $\leq 2n$, cf.\ \cite{Bog,GHJ}. Note that in general $\Sym^nH^2(X,\QQ)(1)\to H^{2n}(X,\QQ)(n)$ is not an isomorphism and even when it is, e.g.\ in the ${\rm K3}^{[2]}$-case, the inclusion $\Sym^nH^2(X,\ZZ)(1)\,\hookrightarrow H^{2n}(X,\ZZ)(n)$ might have a non-trivial finite cokernel, cf.\ \cite{BNWS}.\smallskip

The twist $SH(X,B,\ZZ)$  with respect to a $B$-field lift  $B\in H^2(X,\QQ)$ of $\alpha\in\Br(X)$
is defined as  $SH(X,\ZZ)$ endowed with the pull-back of the Hodge structure
under the composition $$\xymatrix{SH(X,\ZZ)\subset H^{2\ast}(X,\ZZ)\ar[r]^-{\exp(B)}&H^{2\ast}(X,\QQ).}$$ 
It is an easy exercise to check that this does define a Hodge structure, which ultimately
uses $B\in H^2(X,\QQ)\subset SH(X,\QQ)$.
We then define $$\ind_{SH}(\alpha)\coloneqq\ind(SH(X,B,\ZZ)).$$
Note that with this definition, $SH(X,B,\ZZ)$ is a sub-Hodge structure of $H^{2\ast}(X,B,\ZZ)$, the definition of which will be recalled next.

\subsection{Full cohomology}\label{sec:FullCoh}
We next discuss the full even cohomology $H^{2\ast}(X,\ZZ)$, but we will
be brief as topological $K$-theory is in this instance better adapted to our needs.
To obtain a Hodge structure of weight zero we will be actually looking at
$H^{2\ast}(X,\ZZ)(\ast)\coloneqq\bigoplus H^{2k}(X,\ZZ)(k)$ and will ignore possible torsion.\smallskip

 For a fixed $B$-field lift $B\in H^2(X,\QQ)$
of a topologically trivial Brauer class $\alpha\in\Br(X)$ we consider multiplication with
$\exp(B)$ as an isomorphism $\exp(B)\colon H^{2\ast}(X,\QQ)(\ast)\congpf H^{2\ast}(X,\QQ)(\ast)$
of vector spaces
and declare the Hodge structure $H^{2\ast}(X,B,\ZZ)$ to be the one on $H^{2\ast}(X,\ZZ)(\ast)\subset H^{2\ast}(X,\QQ)(\ast)$  obtained by pull-back under $\exp(B)$ of the standard one. In other words, there exists an isomorphism of $\QQ$-Hodge structures of weight zero
$$\exp(B)\colon H^{2\ast}(X,B,\ZZ)\otimes\QQ=H^{2\ast}(X,B,\QQ)\congpf H^{2\ast}(X,\QQ)(\ast).$$
With the standard rank function, i.e.\ the projection onto $H^0(X,\ZZ)$, one then defines
$$\ind_{H^\ast}(\alpha)\coloneqq \ind(H^{2\ast}(X,B,\ZZ)).$$
We will not make much use of this index, but rather use its $K$-theoretic version. To this end,
we next discuss twisting topological $K$-theory.

\subsection{Twisted $K$-theory} 
Assume $X$ that is a smooth complex projective variety and that $\alpha\in\Br(X)$ is a topologically trivial Brauer class,
i.e.\ $\alpha$ is contained in the kernel of the boundary map $\Br(X)\twoheadrightarrow H^3(X,\ZZ)_{\rm tors}$ or,
equivalently, in the image of the exponential map $H^2(X,\QQ)\to\Br(X)$.
 Under this assumption, the topological $K$-theory of $(X,\alpha)$ equals the untwisted topological $K$-theory: $$K_0^{\rm top}(X,\alpha)=K_0^{\rm top}(X).$$
 Via the Chern character, topological $K$-theory is naturally endowed with a Hodge structure of weight zero with the characterising property  that
$$\ch\colon K_0^{\rm top}(X)\otimes\QQ\congpf \bigoplus H^{2i}(X,\QQ(i))$$ is an isomorphism of $\QQ$-Hodge structures of weight zero,  cf.\ \cite{AT}.
Clearly, $K_0(X)\to K_0^{\rm top}(X)$ takes image in the space $K_0^{\rm top}(X)_{\rm Hdg}\subset K_0^{\rm top}(X)$ of integral Hodge classes. Here,
we think of  $K_0(X)$ as the Grothendieck group of $\Db(\Coh(X))$.
\smallskip

Similarly, after fixing a $B$-field lift of $\alpha$, also $K_0^{\rm top}(X,\alpha)$ is endowed with a Hodge structure. More precisely, if $B\in H^2(X,\QQ)$ is a lift of $\alpha\in \Br(X)\subset H^2(X,\ko_X^\ast)
$, then we use the composition
$$\xymatrix{K_0^{\rm top}(X,\alpha)=K_0^{\rm top}(X)\ar[r]^-{\ch}&H^{2\ast}(X,\QQ)\ar[r]^-
{\exp(B)}&H^{2\ast}(X,\QQ)\ar[r]&\bigoplus H^{p,q}(X)}$$
to define a Hodge structure of weight zero on $K_0^{\rm top}(X,\alpha)$ by declaring its
$(i,-i)$-part to be the pre-image of $\bigoplus_{p-q=2i} H^{p,q}(X)$. In other words,
$$\ch^B\coloneqq\exp(B)\cdot \ch\colon K_0^{\rm top}(X,\alpha)\otimes\QQ=K_0^{\rm top}(X)\otimes\QQ\congpf \bigoplus H^{2i}(X,\QQ(i))$$
is an isomorphism of Hodge structures of weight zero.  Note that  choosing another
$B$-field lift of $\alpha$ leads to an isomorphic integral Hodge structure. Hotchkiss \cite[\S\! 3]{HotchHodge}
actually shows
that with this definition $K_0^{\rm top}(X,\alpha)$ is the topological $K$-theory
of the category $\Db(\Coh(X,\alpha))$ with its natural Hodge structure as defined in \cite{Blanc}. 
In principle, one could also work with the twisted 
Hodge structure $H^\ast(X,\alpha,\QQ)=H^\ast(X,B,\QQ)$ as introduced originally in \cite{HS}, but the integral structure of $K_0^{\rm top}(X,\alpha)$ is better adapted to our situation.\smallskip

As in the untwisted case, under the natural map $K_0(X,\alpha)\to K_0^{\rm top}(X,\alpha)$ the Grothendieck group  $K_0(X,\alpha)=K(\Db(\Coh(X,\alpha)))$ is mapped to the space of  integral
Hodge classes
$$K_0^{\rm top}(X,\alpha)_{\rm Hdg}\subset K_0^{\rm top}(X,\alpha).$$ Equivalently, 
$\ch^B(E)\in \bigoplus H^{p,p}(X,\QQ)$ for all $E\in \Db(\Coh(X,\alpha))$, cf.\ \cite[Prop.\ 1.7]{HS}.
 
\subsection{Hotchkiss's Hodge index}\label{sec:HoHo} The rank function is nothing but the composition of the Chern character
with the projection to degree zero,
$$\xymatrix{\rk\colon K^{\rm top}_0(X,\alpha)\ar[r]^-{\ch^B}\ar[r]&\bigoplus H^{2i}(X,\QQ(i))\ar[r]&H^0(X,\QQ)=\QQ.}$$ It takes values in $\ZZ\subset\QQ$. Then, by definition, the index of $\alpha$ is the positive generator
of the image of $\rk$ restricted to the algebraic classes, i.e.
$$\rk(K_0(X,\alpha))=\ind(\alpha)\cdot\ZZ\subset\ZZ.$$
Hotchkiss \cite{HotchHodge} introduced the  Hodge-theoretic index $\ind_H(\alpha)$ by the analogous
condition$$\rk(K_0^{\rm top}(X,\alpha)_{\rm Hdg})=\ind_H(\alpha)\cdot \ZZ\subset \ZZ,$$
where as above $K_0^{\rm top}(X,\alpha)_{\rm Hdg}\subset K_0^{\rm top}(X,\alpha)$ is the group of integral Hodge classes. In other words, $$\ind_H(\alpha)=\ind(K_0^{\rm top}(X,\alpha)).$$The definition only depends on the topologically trivial Brauer class $\alpha$ but not on the choice of the $B$-field lift. For proving lower bounds for the index, see \S\! \ref{sec:Lower}, it is important to observe that
 \begin{equation}\label{eqn:indHind}
 \ind_H(\alpha)\mid\ind(\alpha).
 \end{equation}

The following is essentially due to  \cite{dJP,HotchHodge}, but the result is there not stated
explicitly in this form.

\begin{prop}\label{prop:Shouldbe} Let $X$ be a smooth, complex, projective
variety and assume  $\alpha\in\Br(X)$  is a topologically trivial Brauer class which  can
be represented by a Brauer--Severi variety $\pi\colon P\to X$ for which the integral Hodge conjecture holds.
Then, if $H^\ast(X,\ZZ)$ is torsion free\footnote{Note that for hyperk\"ahler varieties of ${\rm K3}^{[n]}$-type the integral cohomology is known to be torsion free, cf.\ \cite[Thm.\ 1]{MarkmanIntegral}, and no non-trivial torsion has even been detected for any of the other deformation types.} or $\alpha$ is coprime to
$\dim(X)!$, then $$\ind_H(\alpha)=\ind(\alpha).$$
\end{prop}

\begin{proof} Let  $v\in K_0^{\rm top}(X,\alpha)$ be a Hodge class of rank $\rk(v)=\ind_H(\alpha)$. Its pull-back to $P$ gives a Hodge class $\pi^\ast v\in K_0^{\rm top}(P,\pi^\ast\alpha)$.  Now use the direct sum decomposition of algebraic and topological $K$-theory
of the Brauer--Severi variety and its compatibility with pull-back:
$$\xymatrix{K_0^{\rm top}(P)\ar[r]^-\sim& K^{\rm top}_0(P,\pi^\ast\alpha)\ar[r]^-\sim&\bigoplus_{i=0}^{m}K_0^{\rm top}(X,\alpha^i)\\
K_0(P)\ar[u]\ar[r]^-\sim& K_0(P,\pi^\ast\alpha)\ar[u]\ar[r]^-\sim&\bigoplus_{i=0}^{m}K_0(X,\alpha^i),\ar[u]}$$
where we use for both that $\pi^\ast\alpha$ is the trivial Brauer class on $P$.\smallskip

If the integral Hodge conjecture holds for $P$, then  $\pi^\ast v$ is in the image
of $K_0(P)\to K_0^{\rm top}(P)$. Indeed, Perry \cite[Prop.\ 5.16]{Perry}
shows that if $H^\ast(X,\ZZ)$ is torsion free, then the classical integral Hodge conjecture implies the categorical Hodge conjecture. However, then the commutativity of the diagram and the compatibility with the direct sums imply that $v$ is contained in the image of $K_0(X,\alpha)\to K_0^{\rm top}(X,\alpha)$ and hence $\ind(\alpha)\mid\ind_H(\alpha)$.
Combined with the obvious $\ind_H(\alpha)\mid\ind(\alpha)$, this proves the claim.\smallskip

Without the assumption that the cohomology be torsion free, the same argument shows that for
every Hodge class $\gamma\in K_0^{\rm top}(X)$ the class $\dim(X)!^{\dim(X)}\,\gamma$ is algebraic up to torsion.
\end{proof}

\subsection{Comparison}\label{sec:Comp}
For any topologically trivial Brauer class $\alpha\in \Br(X)$ the following Hodge-theoretic
indices have been defined above
$$\ind_{\widetilde H}(\alpha)\,\S\, \ref{sec:twistHM},~\ind_{\Sym^n\widetilde H}(\alpha)\,\S\, \ref{sec:SymHM},~\ind_{SH}(\alpha)\,\S\, \ref{sec:VerbComp},~\ind_{H^\ast}(\alpha)\,\S\, \ref{sec:FullCoh}, \text{ and } \ind_H(\alpha)\,\S\, \ref{sec:HoHo}$$
and we will now clarify the relations between them.\smallskip

\begin{remark}\label{rem:easy}
Here are the easy ones. In the following, $\alpha\in\Br(X)$ will always be a topologically trivial class.\smallskip

$\bullet$ $\ind_H(\alpha)\mid\ind(\alpha)$. This was observed already as (\ref{eqn:indHind}).
In particular, the prime factors of $\ind_H(\alpha)$ all divide $\per(\alpha)$ and, as we will mention below, cf.\ Lemma \ref{lem:LM1}, they actually coincide.\smallskip

$\bullet$ $\ind_{H^\ast}(\alpha)\mid\dim(X)!\cdot \ind_H(\alpha)$. In particular, if $\ind_{H^\ast}(\alpha)$ is coprime to
$\dim(X)!$, then $\ind_{H^\ast}(\alpha)\mid\ind_H(\alpha)$.\footnote{We will later see that the same assertion holds for all $\alpha$ with $\per(\alpha)$ coprime to some integer, possibly larger than $\dim(X)!$.}
 Indeed, $\dim(X)!\cdot\ch\colon K_0^{\rm top}(X)\to H^{2\ast}(X,\QQ)(\ast)$ takes image in the integral cohomology.\smallskip

$\bullet$ $\ind_{\Sym^n\widetilde H}(\alpha)\mid\ind_{\widetilde H}(\alpha)^n$. Indeed, if there exists a Hodge class $v\in \widetilde H(X,B,\ZZ)$ of rank $r$, then its symmetric power $v^{(n)}\in \Sym^n\widetilde H(X,B,\ZZ)$ is a Hodge class of rank $r^n$.\smallskip

$\bullet$ $\ind_{H^\ast}(\alpha)\mid\ind_{SH}(\alpha)$. This follows immediately from $SH(X,B,\ZZ)\subset H^{2\ast}(X,B,\ZZ)$ being an inclusion of Hodge structures.
\end{remark}

Let us first improve upon the second divisibility above.  

\begin{lem}\label{lem:lem1} There exists an integer $N_X$ such that 
$\ind_{H^\ast}(\alpha)=\ind_H(\alpha)$ whenever $\alpha$ is coprime to $N_X$.
\end{lem}

\begin{proof} The inclusion $\dim(X)!\cdot\ch\colon K_0^{\rm top}(X)\,\hookrightarrow H^{2\ast}(X,\ZZ)$
is of finite index, say $k$. Hence, for every Hodge class $v\in H^{2\ast}(X,B,\ZZ)$ the class $k\,v$  is the image
of a Hodge class $w\in K_0^{\rm top}(X,B)$ under $\dim(X)!\cdot\ch$. This shows
$k\cdot\ind_{H^\ast}(\alpha)\mid \dim(X)!\cdot\ind_H(\alpha)$.\smallskip

Now, if $\per(\alpha)$ is coprime to  $N_X\coloneqq k\cdot\dim(X)!$, then also $\ind_H(\alpha)\mid\ind_{H^\ast}(\alpha)$ and hence equality holds.
\end{proof}

Note that so far the various Hodge-theoretic indices fall into two groups: $$\{\ind_H(\alpha),\,\ind_{H^\ast}(\alpha),\,\ind_{SH}(\alpha)\}~\text{ and }~\{\ind_{\widetilde H}(\alpha),\,\ind_{\Sym^n\widetilde H}(\alpha)\}.$$
In order to make contact between the two groups, we need to invoke results of Beckmann \cite{Beck},
Taelman \cite{Taelman}, and Verbitsky \cite{Verb}.  They allow us to view $SH(X,B,\ZZ)$ not only as
a sub-Hodge structure of $H^{2\ast}(X,B,\ZZ)$ but also of $\Sym^n\widetilde H(X,B,\ZZ)$ (up to finite index):
$$H^{2\ast}(X,B,\ZZ)\supset SH(X,B,\ZZ)\subset\Sym^n\widetilde H(X,B,\ZZ).$$
From now on, $n$ will always be fixed by the condition $\dim(X)=2n$.

\begin{lem}\label{lem:lem3}
For all $\alpha\in \Br(X)$ one has $$\ind_{\Sym^n\widetilde H}(\alpha)\mid n!\cdot\ind_{SH}(\alpha).$$
Moreover, there exists a positive integer $N_X$ such that for all Brauer classes coprime to $N_X$ equality holds
$$\ind_{SH}(\alpha)=\ind_{\Sym^n\widetilde H}(\alpha).$$ \end{lem}

\begin{proof} Relying on results of Verbitsky, cf.\ \cite{Bog,Verb}, Taelman \cite{Taelman}  proved that there exists an isometric inclusion of rational Hodge structures
$SH(X,\QQ)\,\hookrightarrow \Sym^n\widetilde H(X,\QQ)$. Moreover, the inclusion
is compatible with the action of the LLV algebra. In fact, by sending $1\in H^0(X,\ZZ)$ (up to scaling) to the symmetric power $(1/n!)\,e^{(n)}$, this compatibility determines the map. Hence, there is an isomorphism of rational Hodge structures 
\begin{equation}\label{eqn:Verb}
SH(X,\QQ)\oplus SH(X,\QQ)^\perp\cong \Sym^n\widetilde H(X,\QQ),
\end{equation}
which is compatible with twisting, i.e.\  for any $B\in H^2(X,\QQ)$, the decomposition
(\ref{eqn:Verb}) can also be read as
\begin{equation}\label{eqn:Verb2}
SH(X,B,\QQ)\oplus SH(X,B,\QQ)^\perp\cong \Sym^n\widetilde H(X,B,\QQ).
\end{equation}

The inclusion $SH(X,B,\QQ)\subset \Sym^n\widetilde H(X,B,\QQ)$ need not 
induce an injection between the underlying integral structures, but there exists a finite index
sublattice $\Gamma\subset SH(X,\ZZ)$, say of index $k$, that is a contained in $\Sym^n\widetilde H(X,\ZZ)$. This implies $\ind(\Sym^n\widetilde H(X,B,\ZZ))\mid \ind(\Gamma)\mid k\cdot \ind(SH)$.
In particular, if $\per(\alpha)$ is coprime to $k$, then $\ind_{\Sym^n\widetilde H}(\alpha)\mid\ind_{SH}(\alpha)$. In fact, using the compatibility with the action of the LLV algebra one finds that we can choose $k=n!$.
Indeed, the Lefschetz operator $e_\lambda$ on $\widetilde H(X,\QQ)$ associated with
$\lambda\in H^2(X,\ZZ)$ is defined over the integers and so is its action on the symmetric power,
see \cite[\S\! 2.1]{Taelman}.\smallskip

For the converse $\ind_{SH}(\alpha)\mid\ind_{\Sym^n\widetilde H}(\alpha)$, we use that $\Gamma\oplus\Gamma^\perp\subset \Sym^n\widetilde H(X,B,\ZZ)$ is a finite index inclusion, say of index $j$. Then for any integral (Hodge) class $v\in \Sym^n\widetilde H(X,B,\ZZ)$ the $\Gamma_\QQ$-component of $j\cdot v$ is contained in $\Gamma$. This then implies $\ind_{SH}(\alpha)\mid j\cdot \ind_{\Sym^n\widetilde H}(\alpha)$ and, therefore, $\ind_{SH}(\alpha)\mid  \ind_{\Sym^n\widetilde H}(\alpha)$ for Brauer classes coprime to $j$. \smallskip

This proves the assertion with $N_X=k\cdot j$.
\end{proof}

We summarise the above discussion by the following.

\begin{prop}\label{prop:summ}
Assume $X$ is a hyperk\"ahler variety of dimension $2n$. Then there exists a positive integer $N_X$
such that for every  $\alpha\in \Br(X)$ coprime to $N_X$ we have
$$\ind_H(\alpha)=\ind_{H^\ast}(\alpha)=\ind_{SH}(\alpha)=\ind_{\Sym^n\widetilde H}(\alpha)\mid\ind_{\widetilde H}(\alpha)^n.$$
\end{prop}

\begin{proof}
We may assume that $|H^3(X,\ZZ)_{\rm tor}|$ divides $N_X$ and then
the assertion only concerns topologically trivial Brauer classes, so that $\ind_{\widetilde H}(\alpha)$
and $\ind_{\Sym^n\widetilde H}(\alpha)$ are well defined.\smallskip

The first and third equalities are the content of Lemma \ref{lem:lem1} and  Lemma \ref{lem:lem3}, while the divisibility assertion has been observed  already in Remark \ref{rem:easy}. There we also observed that $\ind_{H^\ast}(\alpha)\mid\ind_{SH}(\alpha)$. However, similar to the arguments above, the inclusion $SH(X,\ZZ)\oplus SH(X,\ZZ)^\perp\subset H^{2\ast}(X,\ZZ)$ is of finite index, say $k$, and
together with the fact that the rank function is trivial on $SH(X,\ZZ)^\perp$ this shows
that $\ind_{SH}(\alpha)\mid k\cdot \ind_{H^\ast}(\alpha)$, which gives equality for all classes coprime to $k$.\end{proof}

Note that the discussion in this section, especially the last proposition, conclude the proof of Theorem \ref{thm:thm3}.

\section{Period versus Hodge-theoretic indices}

We show that the index is bounded from below for all but possibly one of the Hodge-theoretic indices introduced before. Also, similar to de Jong's result \cite{dJSurf} that period and index coincide for surfaces, we shall prove that period and index often coincide for $\widetilde H$. This will eventually allow us to prove lower bounds for the index in sufficiently many cases.

\subsection{The period divides all indices} The classical fact that $\per(\alpha)\mid\ind(\alpha)$ also holds for the Hodge-theoretic indices. This is the following easy lemma.

\begin{lem}\label{lem:LM1}
Let $\alpha\in \Br(X)$ be a topologically trivial Brauer class on a hyperk\"ahler variety $X$.
Then
 $$\per(\alpha)\mid\ind_{H^\ast}(\alpha)\mid\ind_{SH}(\alpha),~~\text{ } \per(\alpha)\mid\ind_H(\alpha),~\text{ and }~\per(\alpha)\mid\ind_{\widetilde H}(\alpha).$$
 Moreover, if $\alpha$ is coprime, then also $\per(\alpha)\mid\ind_{\Sym^n\widetilde H}(\alpha)$.
\end{lem}

\begin{proof} To prove the first assertion, we fix a Brauer--Severi variety $\pi\colon P\to X$ representing $\alpha$ and a $B$-field lift $B\in H^2(X,\QQ)$. The pull-back map
defines a morphism of Hodge structures $$H^{2\ast}(X,B,\ZZ)\to H^{2\ast}(P,\pi^\ast B,\ZZ)\cong H^{2\ast}(P,\ZZ),$$ where the isomorphism is a consequence of $\pi^\ast \alpha$ being trivial. By construction, this map sends
$v=(r,\ldots)\in H^{2\ast}(X,B,\ZZ)$  to a class $\tilde v\in H^{2\ast}(P,\ZZ)$ which in degree two
is fibrewise $r\,h\in H^2(P_x,\ZZ)\cong\ZZ$. Hence, if $v$ is a Hodge class, then there exists a class in $H^{1,1}(P,\ZZ)$ which fibrewise is of degree $r$. By the Lefschetz $(1,1)$ theorem, this is equivalent to saying that there exists
a line bundle $\ko(r)$ of relative degree $r$. The latter implies $\per(\alpha)\mid r$, see e.g.\ \cite[Lem.\ 1.1]{GH} and, therefore, $\per(\alpha)\mid\ind_{H^\ast}(\alpha)$.\smallskip

The second  assertion $\per(\alpha)\mid\ind_H(\alpha)$ follows,
at least for $\per(\alpha)$ coprime to  $\dim(X)!$, from Lemma  \ref{lem:lem1},  or without any assumption from a result of Hotchkiss \cite[Lem.\ 5.8(4)]{HotchHodge}.\smallskip

To prove $\per(\alpha)\mid\ind_{\widetilde H}(\alpha)$, we first choose a field lift $B$ of $\alpha$ and then observe that
$(r,a,s)\in \widetilde H(X,B,\ZZ)$ is a Hodge class if and only if it is orthogonal to $(0,\sigma,-q(\sigma,B))$.
The latter is in turn equivalent to $a-r\,B\in \NS(X)\otimes\QQ$ and we have to show that in this case
$\per(\alpha)\mid r$.\smallskip

To this end, one observes that $\per(\alpha)$ is the smallest positive integer
$m$ such that the image of $m\,B\in H^2(X,\QQ)$ under the orthogonal projection
$H^2(X,\QQ)\twoheadrightarrow T(X)_\QQ=H^2(X,\QQ)/\NS(X)_\QQ$ is contained in $T'(X)\coloneqq
H^2(X,\ZZ)/\NS(X)\subset T(X)_\QQ$. If $r$ is as above,
then clearly $rB$ has this property and, therefore, $\per(\alpha)\mid r$.\smallskip

The very last assertion is a consequence of Lemma \ref{lem:lem3}.
\end{proof}

\subsection{Period versus index for the Hodge--Mukai lattice} Consider the finite index inclusion
$$\NS(X)\oplus T(X)\subset H^2(X,\ZZ).$$
For any given topologically trivial Brauer class $\alpha\in \Br(X)$ one can choose
a $B$-field lift $B=(1/\ell)\,b\in T(X)\otimes\QQ$ with $b\in T(X)$ primitive and $\ell$ a positive integer.
Then the period of $\alpha$ divides $\ell$ and in many cases, namely in the non-special case, equals it, cf.\ Remark \ref{rem:ellper}.

\begin{lem}\label{lem:inpeH} Assume $X$ is a hyperk\"ahler variety.
Then for any  $\alpha\in\Br(X)$ we have
$$\ind_{\widetilde H}(\alpha)=\per(\alpha).$$
 \end{lem}

\begin{proof} Clearly, $\per(\alpha)\mid\ell$, but in general this is not an equality. Indeed, the period
has the property that for $\ell'\coloneqq \ell/\per(\alpha)$ the class $(1/\ell')\,b\in T(X)_\QQ=T'(X)_\QQ$
is contained in $T'(X)\subset T(X)_\QQ$, where the notation is as above. In other words, there exists
a class $g\in \NS(X)_\QQ$ such that $a\coloneqq (1/\ell')\,b+g\in H^2(X,\ZZ)$. Then $(\per(\alpha),-a,0)\in \widetilde H(X,B,\ZZ)$ is a Hodge class. Hence, $\ind_{\widetilde H}(\alpha)\mid\per(\alpha)$ and the
other direction was already proved in the previous lemma.
\end{proof}

Note that the lemma is consistent with the spirit of the result of de Jong \cite{dJSurf} for surfaces,
which in this light can be seen as a twisted Lefschetz $(1,1)$ theorem, cf.\   Proposition \ref{prop:Shouldbe}.

\begin{remark}\label{rem:ellper}
Recall that a topologically trivial class $\alpha\in\Br(X)$ is called \emph{non-special} if it is associated with a transcendental $B$-field $B=(1/\ell)\,b$
with $\gcd(\ell,q(b))=1$, see \cite{HMSYZ}.  Note that the condition is independent of the choice of the transcendental $B$-field.\smallskip

Furthermore, in this case, $\per(\alpha)=\ell$. Indeed, with the notation as in the last proof, we have $(1/\ell')\,q(b)=q(a,b)\in\ZZ$
and, therefore, $\ell'\mid q(b)$. For a non-special $\alpha$, this immediately gives $\ell'=1$, i.e.\
$\per(\alpha)=\ell$.
\end{remark}

\subsection{Symmetric products}\label{sec:SymProd}
We come back to the divisibility $\ind_{\Sym^n\widetilde H}(\alpha)\mid\ind_{\widetilde H}(\alpha)^n$,
see Remark \ref{rem:easy}. This is not an equality in general, but we will
provide many examples when it is. (The main reason for its failure is the fact that  the index of the inclusion (\ref{eqn:NXB}) below heavily depends on $\per(\alpha)$.)\smallskip

In the following,  $N(X,B)$ denotes the space of all integral Hodge classes in $\widetilde H(X,B,\ZZ)$ and the transcendental lattice $T(X,B)$ is defined as its orthogonal complement. Then the natural inclusion \begin{equation}\label{eqn:NXB}
N(X,B)\oplus T(X,B)\subset \widetilde H(X,B,\ZZ)
\end{equation}
is of finite index.\smallskip

The symmetric product $\Sym^n\widetilde H(X,B,\ZZ)$ is again a 
Hodge structure of weight zero, cf.\ Section \ref{sec:SymHM}, and its lattice of Hodge classes contains $\Sym^nN(X,B)\subset \Sym^n\widetilde H(X,B,\ZZ)$. However, further Hodge class may exist. For example, for $j$ even the quadratic form $q_T$ leads to a Hodge class in $\Sym^jT(X,B)_{\rm Hdg}\otimes\QQ$, which we will
denote by $q^{j/2}_T$. This class has the property that $q^{(j/2)}(\gamma^{(j)},q^{j/2}_T)=q(\gamma)^{j/2}$ for the symmetric power $\gamma^{(j)}\in \Sym^n T(X,B)$ of any class $\gamma\in T(X,B)$.\smallskip

We introduce the following shorthands: $N\coloneqq N(X,B)$, $T\coloneqq T(X,B)$, and $\widetilde H\coloneqq \widetilde H(X,B,\ZZ)$. Then the space of all rational Hodge classes in $H_\QQ^{\otimes n}$ decomposes as
\begin{equation}\label{eqn:SymnT}
\bigoplus_j (N^{\otimes n-j}\otimes\QQ)\otimes (T^{\otimes j}\otimes\QQ)_{\rm Hdg}\cong (H^{\otimes n}\otimes\QQ)_{\rm Hdg}.
\end{equation}
This leads to the finite index inclusion on the level of  the lattices:
$$\bigoplus_j N^{\otimes n-j}\otimes T^{\otimes j}_{\rm Hdg}\subset H^{\otimes n}_{\rm Hdg}.$$

\begin{prop}\label{prop:MTind}
Assume the transcendental lattice $T(X)$ is of rank at least six and
\begin{equation}\label{eqn:assump}
\Sym^jT(X)_{\rm Hdg}\otimes\QQ=\begin{cases}0 &\text{\rm if } j \text{ \rm odd}\\q_T^{j/2}\QQ&\text{\rm if } j\text{ \rm is even.} \end{cases}
\end{equation}
Then for any  non-special coprime Brauer class $\alpha\in\Br(X)$ one has
$$\ind_{\Sym^n\widetilde H}(\alpha)=\per(\alpha)^n.$$
\end{prop}

\begin{proof} Consider a $B$-field lift $B=(1/\ell)\, b\in T(X)$ of $\alpha$, so  $b\in T(X)$ and $\gcd(\ell,q(b))=1$. Then $\ell=\per(\alpha)$, see Remark \ref{rem:ellper}.\smallskip

First observe that  $\tilde e\coloneqq (\ell,-b,0)$ is contained in $N(X,B)$, as it is orthogonal to the generator $(0,\sigma,-q(\sigma,b/\ell))\in T(X,B)^{2,0}$. Hence, the symmetric power $
\tilde e^{(n)}=\ell^ne^{(n)}+\ell^{(n-1)}\cdot b+\cdots+b^{(n)}$ is contained in $\Sym^nN(X,B)\subset  \Sym^n\widetilde H(X,B,\ZZ)_{\rm Hdg}$ and, therefore, $\ind_{\Sym^n\widetilde H}(\alpha)\mid \ell^n$ for all $\ell$ (without any coprimality assumption).\smallskip

We now have to prove that conversely $\ell^n\mid\ind_{\Sym^n\widetilde H}(\alpha)$ for all $\ell$
coprime to a certain $N_X$. To this end, assume that $v\in \Sym^n\widetilde H(X,B,\ZZ)_{\rm Hdg}$ and write
$v=\mu\,\tilde e^{(n)}+v'$ with $\mu\in\QQ$ and $\rk(v')=0$. Then $\mu\,\ell^n=\rk(v)$. It suffices to show that $\mu\in\ZZ$.\smallskip

To illustrate the idea let us first discuss the easier case that $v\in \Sym^n N(X,B)$, in which case we set $N_X=q(b)$. Then 
\begin{equation}\label{eqn:oddeven}
q^{(n)}(b^{(n)},v)=q^{(n)}(b^{(n)},\mu\, \tilde e^{(n)})=\mu\, q(b)^n
\end{equation} and, therefore, $\mu\, q(b)^n\in \ZZ$. The first equality follows from the decomposition, cf.\ the arguments in the proof of Lemma \ref{lem:inpeH},
\begin{equation}\label{eqn:NXB1}
N(X,B)=\ZZ\tilde e\oplus \NS(X)\oplus \ZZ (0,0,1).
\end{equation}
and  the fact that the transcendental class $b$ is orthogonal to $\NS(X)$.
On the other hand, $\mu\,\ell^n=\rk(v)\in\ZZ$ and $\gcd(\ell,q(b))=1$ if $\ell$ is coprime to $N_X=q(b)$. Hence, $\mu\in \ZZ$, which implies
$\ell^n\mid\ind_{\Sym^n\widetilde H}(\alpha)$.
\smallskip

In the general situation, the above argument breaks down in (\ref{eqn:oddeven}),
for $v$ might now have a transcendental part, which eventually forces us to increase $N_X$.
So we decompose the given Hodge class $v$ with respect to (\ref{eqn:SymnT}) as $$v=\sum u_{n-i}\otimes w_i.$$ Note that the classes  $u_{n-i}\in N^{\otimes n-i}\otimes \QQ$ and $w_i\in (T^{\otimes i}\otimes\QQ)_{\rm Hdg}$, are a priori only rational(!) Hodge classes and that $\rk(u_n)=\rk(v)=\mu\,\ell^n\in\ZZ$. Also note that by assumption (\ref{eqn:assump}), $w_i=0$ for $i$ odd and $w_i=\lambda_i\, q_T^{i/2}$ for $i$ even and some scalar $\lambda_i\in \QQ$.
Hence, normalising the $u_{n-i}$, we can write $$v=\sum^{m}_{j=0} (u_{n-2j}\otimes q^{j}_T),$$
where $m$ is determined by $n=2m+1$ resp.\ $n=2m$.\smallskip

We claim that there exists a class $w\in \Sym^n\widetilde H(X,B,\ZZ)$
with the following property
\begin{equation}\label{eqn:Propb}
q^{(n)}(w,v)=\varepsilon\,q^{(n)}\!\left(b^{(n)},u_n\right).
\end{equation}
Here,  $\varepsilon $ is any non-zero integer, which will be fixed below.
 This would immediately lead to $$\varepsilon\,\mu\, q(b)^n=\varepsilon \,q^{(n)}\!\left(b^{(n)},u_n\right)
=q^{(n)}\!\left(w,v\right)\in\ZZ.$$ If now $\ell$ is assumed to be coprime to $N_X\coloneqq 2\varepsilon\,{\rm disc}(T(X))\,q(b)$, then $\mu\,\ell^n\in\ZZ$ and $\varepsilon\,\mu\, q(b)^n\in\ZZ$ together imply
 $\mu\in\ZZ$. \smallskip
 
 It remains to construct $w$. For this we first pick $\eta\in T(X)\cap b^\perp\subset T(X,B)$ with $0<q(\eta)$
 and set $\varepsilon =q(\eta)^m$. Below in Lemma \ref{lem:isotrop} we show that
 under the assumption on the rank of $T(X)$  we can choose $\eta$ with $q(\eta)\mid 8q(b)^2\,{\rm disc}(T(X))^2$.
 Then we define
 \begin{equation}\label{eqn:Defw}
 w\coloneqq \sum_{k=0}^m(-1)^{k}\binom{n}{2k}\,q(\eta)^{m-k}\,q(b)^k\,(b^{(n-2k)}\otimes\eta^{(2k)}).
 \end{equation}
 Setting $\gamma \coloneqq \sqrt{q(\eta)}\,b + i\sqrt{q(b)}\,\eta \in H^2(X,\CC)\subset\widetilde{H}(X,\CC)$, we have that $w = \operatorname{Re}(\gamma^{(n)})$ if $n$ is even, and $w = \frac{1}{\sqrt{q(\eta)}}\operatorname{Re}(\gamma^{(n)})$ if $n$ is odd. 
 The verification of (\ref{eqn:Propb}) is straightforward. 
 First, recall the following formula, which holds for every $u \in \Sym^k\widetilde H,\ v \in \Sym^{n-k}\widetilde H$
 \begin{equation}\label{eq:SplittingPairing}
     q^{(n)}\left(\gamma^{(n)},u \otimes v\right) = \binom{n}{k}\,q^{(k)}\left(\gamma^{(k)},u\right)q^{(n-k)}\left(\gamma^{(n-k)},v\right).
 \end{equation}
 Therefore, 
 \[
    q^{(n)}\left(\gamma^{(n)}, u_{n-2i} \otimes q_T^i\right) = \binom{n}{2i}\,q^{(n-2i)}\left(\gamma^{(n-2i)},u_{n-2i}\right)q^{(2i)}\left(\gamma^{(2i)},q_T^i\right) 
 \]
 which vanishes for $i > 0$ because $q(\gamma) =0 $. 
 It follows that
 \[
    q^{(n)}\left(w,v\right) = q^{(n)}\left(w,u_{n}\right) = q^{(n)}\left(q(\eta)^mb^{(n)},u_{n}\right)
 \]
 where the first equality is obtained by taking the real part from the above, and for the second equality we used that $\eta \in T(X)\cap b^\perp$.

For the last assertion just observe that $\per(\alpha)\mid \ell$ if $\alpha$ has a $B$-field lift of
the form $(1/\ell)\,b$ with $b$ an integral class. 
\end{proof}

\begin{lem}\label{lem:isotrop}
If  $T(X)$ is of rank at least six, then for any $b\in T(X)$ with $q(b)\ne0$ the orthogonal
complement $b^\perp \subset T(X)$ contains a class $\eta$ with $0<q(\eta)\mid 8 q(b)^2 \,{\rm disc}(T(X))^2$. 
\end{lem}

\begin{proof} By Meyer's theorem the indefinite lattice $\Gamma\coloneqq T(X)\cap b^\perp$, which is of rank at least five, contains an isotropic class $0\ne \gamma\in \Gamma$ which we may assume to be primitive. We denote its divisibility by $d$, i.e.\ $q(\gamma,\Gamma)=d\ZZ$ and consider the class
$\eta=2d\,\delta+(2-q(\delta))\,\gamma$, where $\delta\in\Gamma$ is chosen such that
$q(\gamma,\delta)=d\ZZ$. Then $q(\eta)=8d^2$.\smallskip

Then use that $d\mid {\rm disc}(\Gamma)\mid q(b)\,{\rm disc}(T(X))$. For the first divisibility
use that $(1/d)\,\gamma\in \Gamma^\ast$ and also $({\rm disc}(\Gamma)/d)\,\gamma\in \Gamma$.
For the second use the formula ${\rm disc}(\ZZ b\oplus\Gamma)={\rm disc}(T(X))\cdot(T:\ZZ b\oplus \Gamma)$ and the inclusion $T(X)/\ZZ b\oplus \Gamma\,\hookrightarrow (\ZZ b)^\ast/\ZZ b\cong \ZZ/ q(b)\ZZ$.
\end{proof}

%



\section{Upper and lower bounds}\label{sec:Lower}

In this section we prove Theorems \ref{thm:1} and \ref{thm:2}.\smallskip

The upper bound of the Hodge theoretic index provided by Theorem \ref{thm:1} is
an almost immediate consequence of the previous sections, see \S\! \ref{sec:ProofThm1}.
We also comment on further Hodge theoretic obstructions which we currently do not control, see \S\! \ref{sec:furtherobs}.\smallskip

The interest in the  lower bound in Theorem \ref{thm:2} is motivated by the following question: In the theory of hyperk\"ahler varieties, the second cohomology $H^2(X,\ZZ)$ with its K3 Hodge structure
and the quadratic form $q_X$ replaces $H^2(S,\ZZ)$ of a (K3) surface and, indeed, in many respects
hyperk\"ahler varieties often behave more like surfaces than higher-dimensional K\"ahler manifolds. In this spirit one could wonder whether possibly the period-index conjecture advocated for in \cite{HuyPI},
see Conjecture \ref{conj:PIHK}, is  not optimal. Maybe one should rather expect $\per(\alpha)=\ind(\alpha)$ as for surfaces? Theorem \ref{thm:2} shows that this is not the case.

\subsection{Proof of Theorem \ref{thm:1}}\label{sec:ProofThm1} As before, $n$ is fixed such that $\dim(X)=2n$. Then, according to Proposition \ref{prop:summ}, we already know
$$\ind_H(\alpha)=\ind_{H^\ast}(\alpha)=\ind_{SH}(\alpha)=\ind_{\Sym^n\widetilde H}(\alpha)\mid\ind_{\widetilde H}(\alpha)^n$$ for all coprime classes $\alpha$. Combined with Lemma \ref{lem:inpeH},
which shows $\ind_{\widetilde H}(\alpha)=\per(\alpha)$, this proves the claimed 
$\ind_H(\alpha)\mid\per(\alpha)^{\dim(X)/2}$.\smallskip

We conclude the proof by applying Proposition \ref{prop:Shouldbe}.\qed

\subsection{Other Hodge theoretic obstructions}\label{sec:furtherobs}
Theorem \ref{thm:1} says that there is no a priori Hodge theoretic obstruction
to Conjecture \ref{conj:PIHK}. But this is not quite the full picture. Adapting an approach first discussed by  de Jong and Perry \cite{dJP} from the case of general projective varieties to
the hyperk\"ahler situation, there are (at least) two other potential Hodge theoretic obstructions:
Assume $\pi\colon P\to X$ represents a topologically trivial Brauer class $\alpha\in \Br(X)$ such that its index satisfies $\ind(\alpha)\mid\per(\alpha)^{\dim(X)/2}$. This gives rise to two integral Hodge classes:

\begin{enumerate}
\item[(i)]  There exists a closed subvariety
$Z\subset P$ which is generically finite over $X$ of degree $\per(\alpha)^{\dim(X)/2}$.
We denote its cycle class by $\delta\in H^{r,r}(P,\ZZ)$, where $r$ is the relative
dimension of $\pi$. The class $\delta$ has the distinguished property that its restriction to a fibre satisfies $\delta|_{P_x}=\per(\alpha)^{\dim(X)/2}\in H^{2r}(\PP^r,\ZZ)\cong\ZZ$.
\item[(ii)] There exists a closed subvariety $P'\subset P$ of relative dimension $r'\coloneqq\per(\alpha)^{\dim(X)/2}-1$ over $X$ which over a dense open subset of $X$ is linear.
We denote its cycle class by $\Delta\in H^{c,c}(P,\ZZ)$. Here, $c=r-r'$ is the codimension of $P'\subset P$. The class $\Delta$ has the distinguished property that
$\Delta|_{P_x}=h^c\in H^{2c}(\PP^{r'},\ZZ)\cong \ZZ h^c$. (We are taking the dual point of view compared to \cite[Lem.\ 5.13]{dJP}. The class $\gamma$ there corresponds to the dual of $P'$ inside the dual of $P$.)
\end{enumerate}
Geometrically, passing from  (i) to (ii) is realised by taking fibrewise the linear span of the fibre $Z_x\subset P_x$ which yields $P'\subset P$. For the converse direction one uses the zero set of a relative
vector field. Note however, that there is no cohomological way that would describe the transition from
$\delta$ to $\Delta$ or back. Similarly, there is no direct procedure to pass from
the Hodge theoretic index of Hotchkiss to the class $\delta$. See \cite[Lem.\ 5.12]{dJP} for the passage to $\Delta$ (or rather their $\gamma$).\smallskip

In the general situation, so with the exponent $\dim(X)-1$ instead of our hyperk\"ahler exponent $\dim(X)/2$, de Jong and Perry show that there is no obstruction to the existence of the
class $\delta$, see. \cite[Rem.\ 5.22]{dJP}. In contrast, the existence of the class $\Delta$
has consequences for the integral cohomology of $X$ which potentially might 
obstruct the validity of Conjecture \ref{conj:PI}, cf.\ \cite[Rem.\ 5.15]{dJP} and \cite[Thm.\ 1.4]{HotchHodge}.\smallskip

In the hyperk\"ahler situation, the existence of both classes, $\delta$ and $\Delta$, puts  constraints on the cohomology of the hyperk\"ahler variety $X$ which at this point are not known
to hold in general. We briefly explain this in both cases.\smallskip

(i) According to \cite[Prop.\ 5.21]{dJP}, the class $\delta$ exists if and only if there are rational Hodge
classes $c_j\in H^{j,j}(X,\QQ)$, $j=1,\ldots,\dim(X)$,  such  for $i=1,\ldots,\dim(X)$ the classes
\begin{equation}\label{eqn:dJP1}
\ell^{\dim(X)/2-i}\binom{r}{i}\,b^i+\sum_{j=1}^i\binom{r-j}{i-j}\ell^{j-i}\,b^{i-j}c_j
\end{equation}
are integral. Here, $\alpha\in \Br(X)$ is a non-special Brauer class represented by $(1/\ell)\,b\in T(X)_\QQ$ with $\ell=\per(\alpha)$ and $r\geq\dim(X)$ is the relative dimension of a Brauer--Severi
variety $P\to X$ representing $\alpha$.\smallskip

For $1\leq i\leq \dim(X)/2$ the conditions (\ref{eqn:dJP1}) can be easily
satisfied by choosing $c_i=0$. However, this is no longer true for $\dim(X)/2<i$ and in search for a universal solution, the only reasonable choice
for the classes $c_i$ seems to be powers of the BBF-quadratic form $q\in H^{2,2}(X,\QQ)$ up to scaling, i.e.\ $c_i=0$ for $i$ odd
and $c_{2i}=\lambda_i\, q^i\in H^{2i,2i}(X,\QQ)$ for some $\lambda_i\in \QQ$.  
Let us consider the case $\dim(X)=4$ and $i=2$ and $i=3$. The two conditions read
$$\binom{r}{2}b^2+\lambda_1\,q\in H^4(X,\ZZ)~\text{ and }~\frac{1}{\ell}\binom{r}{3}\,b^3+\frac{\lambda_1(r-2)}{\ell}\,(bq)\in H^{6}(X,\ZZ).$$
The first condition roughly amounts to $\lambda_1\in\ZZ$. Concerning the second condition, 
one would like to make use of the fact that $b^3=cq(b)\, (bq)\in H^6(X,\QQ)$ for a certain universal Fujiki constant $c$ over which we have no control in general, so that the condition
comes down to 
$$\frac{1}{\ell}\left(c\binom{r}{3} q(b)+\lambda_1(r-2)\right)\in\ZZ.$$
The existence of an integral solution $\lambda_1$ very much depends on the Fujiki constant $c$.\smallskip

Similarly for the case $\dim(X)=6$ and $i=3$ and $i=4$, one finds the two conditions 
$$\binom{r}{3}b^3+\frac{\lambda_1(r-2)}{\ell}\,(b q)\in H^6(X,\ZZ)~\text{ and }~\frac{1}{\ell}\binom{r}{4}\,b^4+\frac{\lambda_1}{\ell^2}\,(b^2q)+\lambda_2\,q^2\in H^{8}(X,\ZZ),$$
which lead to $\lambda_1\in \frac{\ell}{r-2}\ZZ$ and
$$\frac{1}{\ell^2}\left(\ell c\binom{r}{4} q(b)+\lambda_1\binom{r-2}{2}\right)\in\ZZ.$$
Again, whether the latter equation has a solution $\lambda_1\in\frac{\ell}{r-2}\ZZ$ depend on the  the constant $c$. 
\medskip

(ii) According to \cite[Thm.\ 1.9 \& Lem.\ 5.13]{dJP}, the class $\Delta$ exists if and only if there are rational Hodge
classes $c_j\in H^{j,j}(X,\QQ)$, $j=1,\ldots,\dim(X)$,  such  for $i=1,\ldots,\dim(X)$ the classes
\begin{equation}\label{eqn:dJP2}
\ell^{-i}\binom{\ell^{\dim(X)/2}}{i}\,(-b)^i+\sum_{j=1}^i\binom{\ell^{\dim(X)/2}-j}{i-j}\ell^{j-i}\,(-b)^{i-j}c_j
\end{equation}
are integral, where $\alpha$, $\ell$, and $r$ as above. (These are the conditions for the existence of the class $\gamma$ on the dual Brauer--Severi variety considered in \cite{dJP}, which explains sign of $b$.) Following the same argument as before, this leads to non-trivial conditions on certain Fujiki constants.

\subsection{Mumford--Tate general hyperk\"ahler varieties}\label{sec:MTgen}
The Mumford--Tate group of
a hyperk\"ahler variety $X$ is the smallest subgroup ${\rm MT}(X)\subset {\rm GL}(T(X)\otimes\QQ)$
such that the subgroup ${\mathbb S}^1$ determined by the weight two Hodge structure on $T(X)$
is contained in ${\rm MT}(X)(\RR)$, cf.\ \cite[Ch.\ 3.3.]{HuyK3}.
It has the property that for any $n>0$ the rational Hodge classes in $T(X)^{\otimes n}\otimes\QQ$
are exactly those that are invariant under the induced action of ${\rm MT}(X)$.\smallskip

Assume now that $\kx\to\km$ is the universal family over a component of 
the moduli space of polarised hyperk\"ahler varieties. Then for the very general $t\in \km$,
all fibres $\kx_t$ have the same (maximal) Mumford--Tate group. Such a point $t\in\km$
and the fibre $\kx_t$ are called Mumford--Tate general. According 
to a result of van Geemen and Voisin \cite[Lem.\ 9]{vGV}, the general Mumford--Tate group
is the orthogonal group ${\rm SO}(T(X)_\QQ,q)$. The same result holds as well for the moduli
space of hyperk\"ahler varieties polarised by a lattice of signature $(1,\rho-1)$.\smallskip

As a consequence, one finds that for any Mumford--Tate general (lattice) polarised hyperk\"ahler
variety $X$ the Hodge classes in the tensor algebra $\bigoplus T(X)^{\otimes n}\otimes\QQ$ are exactly
those that are invariant under the action of ${\rm SO}(T(X)\otimes\QQ)$. The same clearly also holds for the symmetric
algebra $\bigoplus \Sym^nT(X)\otimes\QQ$. Classical invariant theory, see e.g.\ \cite[Thm.\ 10.2]{KP}, shows that the
algebra of all symmetric tensors that are invariant under  ${\rm SO}(T(X)\otimes\QQ)$ is the polynomial algebra $\QQ[q_T]$, where $q_T\in \Sym^2T(X)\otimes\QQ$ denotes the class of the quadratic form. Note that this in particular
means that for $j$ odd  $\Sym^jT(X)\otimes\QQ$ does not contain any non-trivial Hodge classes and that for $j$ even it is spanned by $q_T^{j/2}$. Thus, the assumption of Proposition \ref{prop:MTind} are satisfied.

\subsection{Proof of Theorem \ref{thm:2}}\label{sec:PrThmB}
Let first $X$ be an arbitrary hyperk\"ahler variety of dimension $\dim(X)=2n$. Then, by virtue of Proposition \ref{prop:summ},
there exists a positive integer $N_X$ such that 
\begin{equation}\label{eqn:Symn}
\ind_{\Sym^n\widetilde H}(\alpha)=\ind_H(\alpha)
\end{equation}
for all $\alpha\in\Br(X)$ with $\per(\alpha)$ coprime to $N_X$.\smallskip

Let now $X$ be in addition Mumford--Tate general.
According to Proposition \ref{prop:MTind},
\begin{equation}\label{eqn:SymnH}\per(\alpha)^{n}=\ind_{\Sym^n\widetilde H}(\alpha).
\end{equation} for all non-special coprime classes $\alpha\in \Br(X)$. The precise $N_X$ might have to be enlarged at this point.
Now, combining (\ref{eqn:Symn}) and (\ref{eqn:SymnH}) with the obvious $\ind_H(\alpha)\mid\ind(\alpha)$, see (\ref{eqn:indHind}), one concludes $\per(\alpha)^n\mid\ind(\alpha)$.
This concludes the proof of  Theorem \ref{thm:2}.\qed\medskip


\begin{remark}\label{rem:K3n}
Combined with Theorem \ref{thm:thm4}, 
one finds that $\per(\alpha)^{n}=\ind(\alpha)$ holds in this case for
all  non-special coprime classes, see Corollary \ref{cor:equality}.

\end{remark}

\begin{remark}\label{rem:PicardX}
In the proof of Theorem \ref{thm:thm5}, we would like to apply a weaker form of Theorem \ref{thm:2} that only asserts $\per(\alpha)^{\dim(X)/2}\mid\ind(\alpha)$ for sufficiently many classes,
i.e.\ for  classes of arbitrarily large prime period, but
without making any assumption on the Picard number. This is easily achieved by adjusting Proposition  \ref{prop:MTind}. In its proof, the coprime factor $N_X$ depended on $q(\eta)$
with $\eta\in T(X)\cap b^\perp$ which could be bounded only under the assumption that the rank
$T(X)$ is sufficiently large. However, if we simply fix $b$ with $q(b)\ne0$, then we can produce 
Brauer classes of arbitrarily large prime period $p$, namely those associated with $(1/p)\, b$ with $p\nmid q(\eta)\, q(b)$ for an arbitrarily fixed $\eta\in T(X)\cap b^\perp$ of non-vanishing square.
\end{remark}
In the introduction we observed already that the lower bound proved by Theorem \ref{thm:2}
gives even more precise information about the index in all those situations where also an upper bound can be proved. One such a case is described by Theorem  \ref{thm:thm4}, i.e.\ for
hyperk\"ahler varieties of ${\rm K3}^{[n]}$-type, and was recorded
already as Corollary \ref{cor:equality}. Another instance is provided by Lagrangian fibrations due to \cite[Thm.\ 0.3]{HuyPI}.

\begin{cor}\label{cor:LagMT}
Assume $X\to\PP^n$ is a Mumford--Tate general hyperk\"ahler variety with a Lagrangian fibration
such that $\rho(X)\leq b_2-6$.
Then $$\ind(\alpha)=\per(\alpha)^n$$
for all non-special coprime classes $\alpha\in \Br(X)$.\qed
\end{cor}

Recall that the SYZ conjecture predicts that hyperk\"ahler varieties admitting (birationally) a Lagrangian fibration are expected to form a dense countable union of hypersurfaces in the moduli space. These loci can be described as moduli spaces of lattice polarised hyperk\"ahler varieties to
which again \cite[Lem.\ 9]{vGV} applies.



\section{Period-index for ${\rm K3}^{[n]}$-type and prime period}
The goal of this section is to prove Theorem \ref{thm:thm4}. We will prove the
theorem for non-special coprime classes $\alpha$ with $\per(\alpha)=p^m$ a prime
power and for special coprime classes with $\per(\alpha)=p$ a prime number.
This suffices to conclude. Indeed, any Brauer class $\alpha$ with
$\per(\alpha)=\prod p_i^{m_i}$, $p_i\ne p_j$, can be written as 
$\alpha=\alpha_1\cdots\alpha_k$ with $\per(\alpha_i)=p_i^{m_i}$
and then use
$$\xymatrix{\ind(\alpha)\mid \prod\ind(\alpha_i)\mid\prod \per(\alpha)^n.}$$
Here, we tacitly use that with $\alpha$ non-special  also the $\alpha_i$ are non-special,


\subsection{Setup and strategy} The basic idea to prove results about (twisted) bundles  over a hyperk\"ahler variety $X$ of ${\rm K3}^{[n]}$-type is to first consider
the problem for a hyperk\"ahler variety that is actually isomorphic to a Hilbert scheme $S^{[n]}$.
In this case, one can often use the classical K3 surface theory to produce sufficiently many 
(twisted) bundles. Then, as a second step, one tries to establish the same result for the original $X$ by deforming the (twisted or projective) bundles on the Hilbert scheme back to $X$. For the deformation theory
one usually needs the (twisted) bundles on $S^{[n]}$ to be (projectively) hyperholomorphic.\smallskip

There are various construction methods on Hilbert schemes available. The one used by Hotchkiss et al \cite{HMSYZ} is based on work of Markman \cite{Markman}, see also \cite{MSYZ},
and exploits semi-rigid bundles on K3 surfaces. For the majority of cases, however, 
we appeal to results of O'Grady \cite{OG2} which focus on rigid bundles
on K3 surfaces. (Alternatively, one could try to use twisted
bundles obtained by means of Lagrangian fibration as in \cite[\S\! 2]{HuyPI}.) Roughly, the case of special Brauer classes of prime period
will be approached via \cite{HMSYZ,Markman}, see \S\! \ref{sec:special}, while dealing
with non-special Brauer classes  of prime period requires the  existence results proved in \cite{OG2}, see \S\! \ref{sec:nonspecial}.

\subsection{Lattice theory} We start by proving a technical result in lattice theory, which
will be at the heart of the argument in both cases, special and non-special Brauer classes.\smallskip

As we are only dealing with Brauer classes $\alpha\in \Br(X)$ of prime power period, say $p^m$, we
can choose the $B$-field lift of $\alpha$ to be of the form $B=(1/p^m)\,b$ with $b\in T(X)$ not divisible by
$p$. By definition, $\alpha$ is called special if $p\mid q(b)$ and otherwise it is  non-special.
The lattice theory for both cases is similar and covered by the following technical result. The difference between special and non-special is reflected by the choice of $\varepsilon=1,2$ and 
the assertion in (iii).

\begin{lem}\label{lem:lattices}
Let $(X,h)$ be a polarised hyperk\"ahler variety of ${\rm K3}^{[n]}$-type, $p$ a prime number not dividing $q(h)\cdot(2n-2)$, and $b\in T(X)$ a class not divisible by $p$. Then there exist
integers $k\in\{1,\ldots,p-1\}$, $\ell\in\ZZ$ and a class $c\in b^\perp\cap h^\perp$ such that 
$b'\coloneqq \varepsilon k\,b+\ell\, h+\varepsilon p^m\,c$, with $m=1$ if $p\mid q(b)$, satisfies:
\begin{enumerate}
\item[{\rm (i)}] $b'$ is primitive and ${\rm div}(b')=\gcd(\varepsilon,{\rm div}(h))$.
\item[{\rm (ii)}] $q(b')>0$.
\item[{\rm (iii)}] $q(b')\equiv \nu \, q(b)~\text{\rm mod } p$.
\end{enumerate}
Here,  $\nu$ is an arbitrarily fixed element in ${\mathbb F}_p^\ast$ and $\varepsilon=2$ if $p\mid q(b)$ and $\varepsilon=1$ otherwise. In the latter case, one can choose $\ell$ such that $p\nmid\ell$.
\end{lem}

\begin{proof}
Standard lattice theory  shows that there exists a hyperbolic plane $U\subset\langle b,h\rangle^\perp$, cf.\ \cite[Lem.\ 2.2]{HMSYZ} and \cite[Ch.\ 14]{HuyK3}.
Indeed, on the one hand, the saturation $\Gamma\coloneqq\langle b,h\rangle_{\rm sat}\subset H^2(X,\ZZ)$ is a primitive sublattice
of $H^2(X,\ZZ)\cong \Lambda\oplus \ZZ(2n-2)\subset \Lambda\oplus U$, where $\Lambda=E_8(-1)^{\oplus 2}\oplus U^{\oplus 3}$ is the K3 lattice.
On the other hand, there exists a primitive embedding $\Gamma\subset U^{\oplus 2}\subset \Lambda \subset   \Lambda\oplus U$. However, up to automorphisms of $\Lambda\oplus U$ there exists only one primitive embedding $\Gamma\subset \Lambda\oplus U$ and, therefore, we
may assume that $\Gamma$ is contained in $U^{\oplus 2}\subset E_8(-1)^{\oplus 2}\oplus U^{\oplus 2}\subset \Lambda\subset H^2(X,\ZZ)$. In particular, there exists a hyperbolic plane
in $\Gamma^\perp=\langle b,h\rangle^\perp$.\smallskip

Let $e,f\in U$ be the standard isotropic base with $(e,f)=-1$ as before and let $c=e+\lambda f\in U\subset \langle b,h\rangle^\perp$ with  $\lambda$ to be specified later. Note that then for any choice of integers $k$, $\ell$, and $\lambda\gg0$ the class $b'=\varepsilon k\,b+\ell\, h+\varepsilon p^m\,c$ satisfies $q(b')\gg0$  and  ${\rm div}(b')\mid \varepsilon p$. For the latter use $q(b',f)=q(\varepsilon p\,c,f)=-\varepsilon p$. \smallskip

Now choose a class $a\in H^2(X,\ZZ)$ with $q(b,a)={\rm div}(b)$. Note that $p\not\mid {\rm div}(b)$, as otherwise $(1/p)\,b\in H^2(X,\ZZ)^\ast$ and hence, since ${\rm disc}(H^2(X,\ZZ))=(2n-2)$ and $p\not\mid (2n-2)$,  also $(1/p)\,b\in H^2(X,\ZZ)$, which is excluded by assumption.\smallskip

Assume first that $p\mid q(b)$ and so $\varepsilon=2$ and $m=1$. Then we let $\ell=p$ and find $q(b',a)\equiv 2k\, q(b,a)\equiv
2k\, {\rm div}(b)\not\equiv 0\text{ mod }p$ for any $k\in \{1,\ldots,p-1\}$. Hence, ${\rm div}(b')=\gcd(2,{\rm div}(h))$ and, at the same time, $q(b')\equiv q(b)\equiv 0\text{ mod }p$. Obviously, for $p\ne 2$ the vector $b'$ is primitive.\smallskip

Let now $p\not\mid q(b)$ and so $\varepsilon=1$ and $m$ an arbitrary positive interger. Then, for any choice of $\ell$ and $k\in\{1,\ldots, p-1\}$, one has $p\not\mid k\, q(b)=q(b',b)$ and, therefore, ${\rm div}(b')=1$. It remains to verify that $k$ and $\ell$ can be chosen such that  $q(b')\equiv k^2q(b)+\ell^2q(h)\equiv \nu\,q(b)\text{ mod }p$.
As $p\not\mid q(b)\,q(h)$, this follows from the fact that any non-degenerate conic over ${\mathbb F}_p$ admits an ${\mathbb F}_p$-rational point with non-trivial first coordinate.
\end{proof}

\subsection{Proof of Theorem \ref{thm:thm4} for special Brauer classes}\label{sec:special}
The proof in this case is entirely reduced to \cite{HMSYZ} with a minor modification exploiting $p\mid q(b')$, which is ensured by  Lemma \ref{lem:lattices}, (iii). 
We assume that $p\nmid q(h)\cdot n\cdot(2n-2)$. There are two cases to be distinguished:
 ${\rm div}(b')=1$ and ${\rm div}(b')=2$, cf.\ Lemma \ref{lem:lattices}, (i).\smallskip

$\bullet$ Let us first assume ${\rm div}(b')=2$. Following the discussion of  Case (A.2) in the proof of \cite[Prop.\ 2.3]{HMSYZ}, one constructs a parallel transport operator
$$\rho\colon H^2(S^{[n]},\ZZ)\congpf H^2(X,\ZZ), ~-2\, H+p\delta\mapsto \pm b'.$$
Note that the class $L_{u_1}$ there corresponds here to our class $b'$.\smallskip

The second step applies \cite[Prop.\ 1.2 \& 1.3]{HMSYZ}, which we will not recall in full. However, unlike the choice of $r=4p^2$ and $m=\mp 4p$ there we can choose $r=p$ and $m=\mp 1$. The crucial assumption $2r\mid m^2 H^2$ is still satisfied for this choice, due to our assumption that $p\mid q(b')$. Indeed, by definition of $\rho$ we have $4H^2=q(b')+p^2(2n-2)$ and so
$p\mid H^2$, but for $p\ne 2$ this implies $2p\mid H^2$.\smallskip

 As $(1/2p)\,b'\in H^2(X,\QQ)$ induces the Brauer class $\alpha^k$, the crucial \cite[Prop.\ 1.2 \& 1.3]{HMSYZ} implies that $\ind(\alpha^k)\mid n!\cdot p^n$. To conclude, apply \cite[Thm.\ 1.4 (2)]{AW} (recall that $p\nmid k$) and the assumption $p\nmid n$, which together eventually show $\ind(\alpha)=\ind(\alpha^k)\mid p^n$.\smallskip
 
$\bullet$ The case ${\rm div}(b')=1$ is similar. It is covered by the discussion of Case (A.1) 
in  \cite[Prop.\ 2.3]{HMSYZ}, except that now the parallel transport operator is constructed to satisfy
$$\rho\colon H^2(S^{[n]},\ZZ)\congpf H^2(X,\ZZ), ~ H-p\delta\mapsto \pm b'.$$

 This time we choose $r=4p$ and $m=\mp2$ (while $r=4p^2$ and $m=\mp2p$ in
 \cite[Prop.\ 2.3]{HMSYZ}), for which $p\mid q(b')$ again implies $2r\mid m^2H^2$.
The rest of the argument is identical to the above.
\qed

\subsection{Non-special Brauer classes}\label{sec:nonspecial}
We are now turning to the case of non-special Brauer classes of prime power period $p^m$. In other words, to classes  $\alpha\in\Br(X)$ with a $B$-field lift
of the form $B=(1/p^m)\,b$ with $b\in T(X)$ and such that $p\nmid q(b)$. In this case, instead of applying \cite{HMSYZ} which relies on the work of Markman, we use results of O'Grady \cite{OG2}.\smallskip

More precisely, we appeal to \cite[Thm.\ 1.1]{OG2}. We will not recall the full statement here but
to facilitate the translation we record that the numerical invariants 
$r_0, g,l,e,\bar e$  in  \cite[Thm.\ 1.1]{OG2}, will in our situation take the following values: $r_0=p^m$, 
$g=l=1$, and $\bar e=4e=4q(b')$. Here, $b'\in H^2(X,\ZZ)$ is the class provided by Lemma \ref{lem:lattices}, where $\nu\in {\mathbb F}_p^\ast$ will be specified in the proof below.
There is an additional integer $i$ in the assertion of the result, which is  $i=1$, as $p\ne2$.

\begin{lem}
Assume $p\nmid q(h)\cdot(2n-2)\cdot (n+3)$, then $\nu$ can be chosen such that 
with the above values of $r_0, g,l,e$, and $\bar e$  assumptions {\rm (1.2.2)-(1.2.5)}  in \cite[Thm.\ 1.1]{OG2} are satisfied.
\end{lem}

\begin{proof}
The assumption (1.2.2) says $g\mid (r_0-1)/2$ for $r_0$ odd, which clearly holds for $r_0=p^m\ne2$
and $g=1$.  The assumption (1.2.3) comprises the three statements:
$l\mid (n-1)$, $\gcd(l,r_0)=1$, and $\gcd(l,(r_0-1)/g)=1$, which are all trivially satisfied for our choice $l=1$. The assumption (1.2.5) requires $\bar e+2(n-1)(r_0-1)^2/g^2\equiv 0\text{ mod } 8l^2$,
which with our choices amounts to $8\mid 4q(b')+2(n-1)(p^m-1)^2$, a consequence of $q(b')$ and $p^m-1$ both being even.\smallskip

The only one of the assumptions in O'Grady's theorem that needs a moment thought is (1.2.4)
which reads $g^2\bar e+2(n-1)(r_0-1)^2+8\equiv 0\text{ mod } 8r_0$
which becomes $4q(b')+2(n-1)(p^m-1)^2+8\equiv 0\text{ mod } 8p^m$ or, equivalently,
$2q(b')+(n-1)+4\equiv 0\text{ mod } p^m$. According to Lemma \ref{lem:lattices},
$b'=k\, b+\ell\, h+p^m\,c$ can be chosen such that it satisfies this equation $\text{ mod }p$ if $\nu$ is fixed with
$2\nu \,q(b)\equiv -(n+3)\text{ mod }p$. The latter is possible as we are assuming $p\nmid n+3$.\smallskip

In order to lift this equation $\text{ mod } p^m$, we modify $b'$ to $b'+j\,h$ for some appropriate $j$. Note that changing $b'$
in this way does not affect the verifications done before. Using $q(b'+j\,h)=q(b')+q(h)\,j^2+(2\ell q(h))\,j$, the existence of
$j$ such that $2q(b'+j\, h)+(n+3)\equiv 0\text{ mod }p^m$ follows from
Hensel's lemma applied to the quadratic equation $2q(h)\,x^2+(4\ell q(h))\, x+2q(b')+(n+3)$ while using the assumption $p\nmid 2q(h)$ and $p\nmid\ell$ according to Lemma \ref{lem:lattices}.
\end{proof}
\smallskip

\noindent{\bf Proof of Theorem \ref{thm:thm4} for non-special Brauer classes.} Once the assumptions of \cite[Thm.\ 1.1]{OG2} are verified, O'Grady's result asserts the existence
of a $\mu_L$-stable vector bundle $E$ on the general polarised hyperk\"ahler
variety of ${\rm K3}^{[n]}$-type $(Y,L)$ with the following numerical invariants:
$$\rk(E)=(p^m)^n,~{\rm c}_1(E)=(p^m)^{n-1} \,{\rm c}_1(L),\text{ and } \Delta(E)=\frac{(p^m)^{2n-2}\cdot((p^m)^2-1)}{12}\,{\rm c}_2(Y)$$ and such that $H^\ast(Y,\kend_0(E))=0$. 
\smallskip

The polarised $(Y,L)$ is viewed as a deformation of $X$ together with the class $b'\in H^2(X,\ZZ)$.
In particular, the projectivisation $\PP(E)$ of the bundle $E$ deforms to a Brauer--Severi variety $P\to X$ of relative dimension $p^{mn}-1$.
Here we are using that with the given discriminant, which is of type $(2,2)$ on every deformation, the projectivisation of the  stable bundle is hyperholomorphic, see \cite{VerbJAG}.  Recall  that the image of $\PP(E)$ under the natural boundary map $H^1(Y,{\rm PGL}(p^{mn}))\to H^2(Y,\mu_{p^{mn}})$ is the image of $-{\rm c}_1(E)=-p^{m(n-1)}\,{\rm c}_1(L)$ under 
the projection $H^2(Y,\ZZ)\to H^2(Y,\mu_{p^{mn}})$, see e.g.\ \cite[Lem.\ 1.5]{HS}. Thus,
$$(1/p^{mn})\,(p^{m(n-1)}\,b')=(1/p^m)\,b'=(k/p^m)\,b+(\ell/p^m)\,h+c\in H^2(X,\QQ),$$
with   $(\ell/p^m)\,h\in H^{1,1}(X,\QQ)$ and $c\in H^2(X,\ZZ)$, lifts
$\alpha^k$ and corresponds to $(1/\rk(E))\,{\rm c}_1(E)=(1/p^{mn})\,(p^{m(n-1)}L)=(1/p^m)\,L$ under deformation.\smallskip

 Hence, the class of $[P\to X]$ in $H^2(X,\mu_{p^{mn}})\cong H^2(Y,\mu_{p^{mn}})$ is
a lift of $\alpha^k$ under the map $H^2(X,\mu_{p^{mn}})\to \Br(X)[p^{mn}]$. This proves  $\ind(\alpha^k)\mid p^{mn}$, which according to \cite[Thm.\ 1.4 (2)]{AW} also implies $\ind(\alpha)\mid p^{mn}$, as 
$k$ is coprime to $p$.\qed

\section{Covering families of curves}
Due to the work of many people, one knows that every K3 surface $S$ is covered by (singular) elliptic curves, see
\cite[Ch.\ 13.2]{HuyK3} for details and references. More precisely, there exists a curve $B\subset|H|$ in some linear system generically parametrising integral curves of geometric genus one. Equivalently,
there exists a smooth elliptic fibration $\kc\to B$ over an open curve $B$ together with a dominant morphism
$\pi\colon\kc\twoheadrightarrow S$. As explained in \cite{HM}, this leads to an alternative proof of the
period-index conjecture for K3 surfaces, at least for coprime classes. We briefly recall the main steps of the argument.\smallskip

As a first step, by base change to a appropriate finite cover of $B$, we may assume that $\kc\to B$ comes with a section $\widetilde B\subset \kc$. Its image $\pi(\widetilde B)\subset S$ is either a point or a curve, but in any case $\Br(S)\to \Br(\pi(\widetilde B))$ is trivial. Therefore, the pull-back of $\alpha$ to $\kc$, which for simplicity we call again $\alpha\in\Br(\kc)$, has the distinguished property that it restricts trivially to all fibres $\kc_b$ and
to the section $\widetilde B\subset\kc$:
$$\alpha|_{\kc_b}=1~\text{ and }  ~\alpha|_{\widetilde B}=1.$$

Under these assumptions, the moduli space
$\Pic^0_\alpha(\kc/B)$ of $\alpha$-twisted sheaves of rank one on the fibres $\kc_b$ is a fine moduli
space, so there exists a relative $(1,\alpha)$-twisted Poincar\'e sheaf $\kp$ on $\Pic^0_\alpha(\kc/B)\times_B\kc$, see \cite{HM} for a detailed discussion. The induced Fourier--Mukai transform $\Phi\colon\Db(\Pic^0_\alpha(\kc/B))\to \Db(\kc)$ can be used to produce $\alpha$-twisted sheaves of low rank on $\kc$ and we currently have two different ways to do this.\smallskip

Either, one uses the morphism $\varphi\colon \Pic^0_\alpha(\kc/B)\twoheadrightarrow \Pic^0(\kc/B)$, $L\mapsto L^{\per(\alpha)}$, which is \'etale of degree $\per(\alpha)^2$, and the structure sheaf $\ko_Z$
of the pre-image $Z=\varphi^{-1}(0)$ of the zero-section to produce the $\alpha$-twisted sheaf $\kp|_{Z\times_B\kc}=\Phi(\ko_Z)$
of rank $\per(\alpha)^2$. Or, one glues the power $\Theta^{\per(\alpha)}\in \Pic^0(\kc/B)$ of the principal polarisation on the relative Jacobian to a
line bundle $\Theta_\alpha\in \Pic^0_\alpha(\kc/B)$ and considers $\Phi(\Theta_\alpha)$, which is 
an $\alpha$-twisted sheaf of rank $\per(\alpha)$, cf.\ the discussion in  \cite[\S\! 3.2]{HuyPI}.
The second construction method gives a better result. Namely, the index of $\alpha\in\Br(\kc)$
equals the period of $\alpha$. Pushed forward to $S$ one finds that $\per(\alpha)=\ind(\alpha)$ for all $\alpha\in\Br(S)$ coprime to the degree of the covering map $\kc\to S$. 

\subsection{Covering families of curves} In \cite{HM} it was explained how the above considerations can be generalised to the case of arbitrary smooth projective varieties. 
Fixing a very ample linear system on a variety  $X$ of dimension $d$
allows one to construct a family $\kc\to B=\PP^{d-1}$ of complete intersection curves
dominating $X$. Once again, because all $\kc_b\subset X$ are curves, Brauer classes
on $X$ are trivial on the fibres of $\kc\to B$. Furthermore, the family can be constructed in such
a way
that all curves pass through a fixed point $x\in X$, which gives rise to a section $\widetilde B\subset\kc\to B$ to which all Brauer classes on $X$ pull-back trivially. This is all one needs to carry out the above program. The two construction techniques 
sketched above produce $\alpha$-twisted sheaves of rank $\per(\alpha)^{2 g(\kc_b)}$
resp.\ $\per(\alpha)^{g(\kc_b)}$, see \cite[Thm.\ 1.2]{HM} and \cite[\S\! 3.2]{HuyPI},
which eventually leads to $\ind(\alpha)\mid\per(\alpha)^{g(\kc_b)}$ for all classes $\alpha\in\Br(X)$.\smallskip

The fact that the curves $\kc_b$ are complete intersection curves was used only at two places:
First, to produce a family $\kc\to B$ with $\kc$ birational to $X$ and, second, to ensure that all curves
pass through a fixed point $x\in X$. If the first of the two properties is not guaranteed, then
one is in a situation similar to the one for covering families of elliptic curves on K3 surfaces. It has the effect that the upper bound is only proved for classes that are coprime to the degree of
the generically finite map $\kc\to X$. In contrast, the second property is crucial. If the curves do not pass
through a fixed point, we cannot ensure that $\kc\to B$ has a (multi-)section to which 
all (coprime) classes $\alpha\in \Br(X)$ restrict trivially. As a consequence, the Poincar\'e sheaf $\kp$
on $\Pic_\alpha^0(\kc/B)\times_B\kc$ only exists as a $(\beta,\alpha)$-twisted sheaf for some
Brauer class $\beta\in\Br(\Pic_\alpha^0(\kc/B))$ and 
exploiting the Fourier--Mukai functor $\Phi\colon \Db(\Pic_\alpha^0(\kc/B),\beta^{-1})\to\Db(\kc,\alpha)$ would require to first produce suitable $\beta^{-1}$-twisted sheaves on $\Pic_\alpha^0(\kc/B)$, which is as complicated as the original period-index problem for $\kc$ or for $X$.\smallskip

Theorem \ref{thm:thm5}, which will be proved further below, turns this observation into a non-existence result.

\subsection{Elliptic curves on Hilbert schemes}\label{sec:curvesHilb} The Hilbert scheme
$S^{[n]}$ of a K3 surface with a genus one fibration $f\colon S\to\PP^1$ comes itself
with a fibration $S^{[n]}\to\PP^n=(\PP^1)^{(n)}$, the generic fibre of which is a product of elliptic curves $S_{t_1}\times \cdots\times S_{t_n}$, $(t_1,\ldots, t_n)\in (\PP^1)^{n}$. If one fixes $n-1$ generic points $x_1,\ldots, x_{n-1}\in S$ and $t\in \PP^1$, then the locus $\{\,[Z]\mid x_1,\ldots, x_{n-1}\in Z,~t\in f(Z)\,\}\subset S^{[n]}$ describes a curve of genus one isomorphic to the fibre
$f^{-1}(t)$. This construction gives rise to a dominating family
\begin{equation}\label{eqn:coveringell}
\xymatrix@R=15pt{\kc\ar[d]\ar@{->>}[r]& S^{[n]}\\
   U&}
   \end{equation}
   of genus one curves  parametrised by an open dense subset $U\subset S^{n-1}\times \PP^1$.\smallskip
   
The discussion can be generalised to arbitrary K3 surfaces. Instead of a genus one fibration
$S\to \PP^1$
one uses a covering family of (singular) elliptic curves $\xymatrix@C=17pt{B&\ar[l]\kc\ar[r]& S}$
to produce a dominating family (\ref{eqn:coveringell}) of curves of genus one contained in the Hilbert scheme $S^{[n]}$.\smallskip

The Brauer group of $S^{[n]}$ is naturally isomorphic to the Brauer group of $S$. There is a Hodge theoretic argument for this observation as well as a geometric explanation, cf.\ \cite[\S\! 3.4]{HuyPI}. According to  \cite[Thm.\ 0.4]{HuyPI}, this allows one to conclude that $\ind(\alpha)\mid\per(\alpha)^{n}$ for
all $\alpha\in \Br(S^{[n]})=\Br(S)$. Now, on the one hand, the covering family (\ref{eqn:coveringell}) of genus one curves suggests that possibly $\per(\alpha)=\ind(\alpha)$, at least for coprime classes. On the other hand, by Remark \ref{rem:K3n}, we know that  for the Hilbert scheme of a Mumford--Tate general K3 surface of Picard number one (at most 16 would be enough) equality $\ind(\alpha)=\per(\alpha)^n$ holds
all  non-special coprime classes. Theorem \ref{thm:thm5} resolves this `contradiction'. Indeed, the elliptic curves in $S^{[n]}$ constructed in the above way clearly
do not pass through one fixed point.

\subsection{Proof of Theorem \ref{thm:thm5}}\label{sec:proofthm5}
The idea of the proof follows the reasoning above.\smallskip

 On the one hand, for a Mumford--Tate general hyperk\"ahler variety $X$ of Picard number $\rho(X)\leq b_2(X)-6$, Theorem \ref{thm:2} ensures that there does not exist an integer $N_X$
such that for all Brauer classes $\alpha\in\Br(X)$ coprime to $N_X$ the index $\ind(\alpha)$ is strictly smaller than $\per(\alpha)^{\dim(X)/2}$.
\smallskip

 On the other hand, if $\xymatrix{B&\ar[l]\ar[r]\kc\ar[r]&X}$ is a dominating family of smooth curves of genus $g=g(\kc_b)$ all passing through a fixed point $x\in X$, then
 running the program recalled above for families of curves of genus $g$ would yield
 $\ind(\alpha)\mid\per(\alpha)^{g}$ for all coprime classes $\alpha\in\Br(X)$.\smallskip
 
 Both statements together immediately show $g(\kc_b)\geq\dim(X)/2$.\smallskip
 
Without  the assumption on the Picard number of $X$, one uses Remark \ref{rem:PicardX}, which shows that 
 there still exist Brauer classes of arbitrary large prime period satisfying $\per(\alpha)^{\dim(X)/2}\mid\ind(\alpha)$. This is enough for the above argument.\qed

\begin{remark}
It is straightforward to apply the same arguments to rule out other dominating families.
For example, if $\{Y_t\}$ is a covering family of isotropic subvarieties $Y_t\subset X$ all
containing a fixed point $x\in X$, then  every topologically trivial Brauer class on $X$
is contained in the kernel 
of $\Br(X)\to\Br(Y_t)$ and the techniques in \cite{HuyPI,HM} allow one to show that
$\ind(\alpha)\mid\per(\alpha)^{q(Y_t)}$ for all coprime classes $\alpha\in\Br(X)$. Here, $q(Y_t)=h^{0,1}(Y_t)$ is the dimension of the Picard variety of $Y_t$.
 Hence, according to Theorem \ref{thm:2},  $X$ cannot be Mumford--Tate general or $$\dim(X)/2\leq q(Y_t).$$
 
Also note that the condition that the curves $\kc_b$ or the isotropic subvarieties $Y_t$ pass through
one fixed point $x\in X$ can be relaxed. Indeed, it suffices to assume that there exists
a subvariety $Z\subset X$ for which the restriction map $\Br(X)\to \Br(Z)$ is essentially trivial
and such that $Z$ intersects all $\kc_b$ resp.\ $Y_t$. Typically, such $Z$ are produced as isotropic subvarieties.
\end{remark}


\end{document}